\documentclass[final,1p]{elsarticle_}
\usepackage[english]{babel}
\usepackage{amsmath,amssymb,graphicx,color}
\usepackage{mathrsfs}
\usepackage[colorlinks, pdfauthor={Yu. Brezhnev},
            pdfwindowui=false,%pdfmenubar=false,%pdftoolbar=false,
            hyperfootnotes=true, %pagebackref=true,
            bookmarksopen=false,bookmarks=false]{hyperref}

\usepackage{CutPage}
\textcenter{0mm}{0mm}
\croppage{2mm}{-4mm}
\usepackage{YURA,myScripts}
\usepackage{mySymbols}

\newtheorem{theorem}{Theorem}%[section]
\newtheorem{lemma}[theorem]{Lemma}
\newtheorem{corollary}[theorem]{Corollary}
\newtheorem{proposition}[theorem]{Proposition}
\newdefinition{definition}{Definition}
\newdefinition{remark}{Remark}
\newdefinition{example}{\sf Example}
\newproof{pf}{Proof}

\makeatletter % \dddot плохо работает в elsarticle: \larger[3] -> [2]
\def\dddot#1{{\mathop{{}#1}\limits^{\vbox to-1.4\ex@{\kern-\tw@\ex@
	\hbox{\larger[2]$\m@th\mskip0.2mu.\mskip-2.2mu.\mskip-2.2mu.$}%
	\vss}}}}
\makeatother

\def\ded{\oper\eta}

\def\Aleph{{\bo{\aleph}}}
\def\G{{\bo{\mathfrak{G}}}}
\def\bcQ{{\bo{\mathcal{Q}}}}
\def\R{{\mathscr{R}}}
\def\p{\mathcal{P}}
\def\Psix{\mbox{$\mathcal{P}_6$}}
\def\HG{\mathbb{H}^+\!/\G_x}
\def\T{\hbox{\smaller[2]$\mathrm{T}$}}
\newcommand{\SUM}[3]{{\sideset{}{_{#1}}\sum_{#2}^{#3}}}
\journal{Journal of Differential Equations}

\begin{document}
\begin{frontmatter}

\hfill
\href{http://arXiv.org/abs/1011.1645}%
{\ttt{\smaller[2]http://arXiv.org/abs/1011.1645}}
\bigskip

\title{The sixth Painlev\'e transcendent\\
and uniformization of algebraic curves}

\tnotetext[t1]{Research was supported by the Tomsk State University
Competitiveness Improvement Program}

\author{{\sc Yurii~V.~Brezhnev}}
\ead{brezhnev@mail.ru}

\address{Department of Quantum Field Theory\\%
Tomsk State University\\$634050$ Tomsk, Russia}

\begin{abstract}
We exhibit a remarkable connection between sixth equation of
Painlev\'e list and infinite families of explicitly uniformizable
algebraic curves. Fuchsian equations, congruences for group
transformations, differential calculus of functions and
differentials on corresponding Riemann surfaces, Abelian integrals,
analytic connections (generalizations of Chazy's equations), and
other attributes of uniformization can be obtained for these curves.
As byproducts of the theory, we establish relations between
Picard--Hitchin's curves, hyperelliptic curves, punctured tori,
Heun's equations, and the famous differential equation which Ap\'ery
used to prove the irrationality of Riemann's $\zeta(3)$.
\end{abstract}


\end{frontmatter}

\newpage
\tableofcontents

\section{Introduction}

The first example of  general solution to the famous sixth
Painlev\'e transcendent
\begin{equation}\label{p6}
\begin{aligned}
\Psix\colon\quad y_{\mathstrut\mathit{xx}}^{\mathstrut}&=\frac12\!
\left(\frac1y+\frac{1}{y-1}+ \frac{1}{y-x}\right)\!
y_{\mathstrut x}^{2}-\left(\frac1x+\frac{1}{x-1}+
\frac{1}{y-x}\right)\! y_{\mathstrut x}^{\mathstrut}\\
&\== +\frac{y\,(y-1)(y-x)}{x^2(x-1)^2}\!
\left\{\alpha-\beta\,\frac{x}{y^2}+\gamma\,\frac{x-1}{(y-1)^2}-
\mbig[7](\delta-\Mfrac12\mbig[7])\,\frac{x\,(x-1)}{(y-x)^2} \right\}
\end{aligned}
\end{equation}
was obtained in 1889 before equation \eqref{p6} itself had been
derived by Richard Fuchs in 1905 \cite{fuchs}. This case corresponds
to parameters $\alpha=\beta=\gamma=\delta=0$ and is referred
frequently to as Picard's solution \cite{picard}. Surprisingly, but
the second one was obtained by N.~Hitchin \cite{hitchin} after more
than one hundred years. It corresponds to parameters
$\alpha=\beta=\gamma=\delta=\frac18$. Presently these solutions are
the only instances, up to automorphisms in the space
$(\alpha,\beta,\gamma,\delta)$, when solution of \eqref{p6} is known
in its full generality.

One year after Fuchs, Painlev\'e \cite{painleve} gave a remarkable
form to  \eqref{p6} which  is known nowadays as the $\wp$-form of
$\Psix$-equation. The modern representation of this result is given
by the nice equation obtained independently by Babich \& Bordag
\cite{babich} and Manin \cite{manin}:
\begin{equation}\label{p6wp}
-\frac{\pi^2}{4}\,\frac{d^2z}{d\tau^2}=
\alpha\,\wp'(z|\tau)+\beta\,\wp'(z-1|\tau)+\gamma\,\wp'(z-\tau|\tau)+
\delta\,\wp'(z-1-\tau|\tau)\,.
\end{equation}

\subsection{Motivation}

It is not a matter of common knowledge that original motivation of
Picard and Painlev\'e, when deriving their results, was closely
related to construction of single-valued analytic functions.
Painlev\'e himself repeatedly wrote  (see his \OE{}uvres \cite{PO})
about representing solutions in terms of single-valued functions and
Picard mentions it about throughout his almost 200-page-long
treatise \cite{picard}. See also the survey article by R.~Conte in
\cite[p.~77--180]{conte} and book \cite{CM} wherein
single-valuedness of functions are constantly emphasized. It is also
known that original statement of the problem on fixed critical
singularities in solutions of ordinary differential equations
(\odes) was initiated by Picard himself \cite[Ch.~V]{picard} and
subsequently developed by Painlev\'e and Gambier \cite{gambier}.

\subsubsection{Equation \Psix, algebraic solutions, and Riemann
surfaces}

On the other hand, it has long been known that equation \Psix\ is a
rich source  of algebraic solutions/functions
\begin{equation}\label{F}
F(x,y)=0
\end{equation}
and genera of corresponding Riemann surfaces $\R$ may be as great as
is wished. Effective description of such surfaces and single-valued
objects on them is the subject of the theory of uniformization of
algebraic curves and automorphic functions
\cite{lehner,fricke,bers,ford}. In this theory, the algebraic
irrationality \eqref{F} is completely determined by the two
principal transcendental meromorphic objects on $\R$. These are
function field generators $x(\tau)$ and $y(\tau)$, of which all
other meromorphic functions  are built: rational functions $R(x,y)$,
differentials, Abelian differentials, and their integrals.
Fundamentally, the constructive `meromorphic analysis' on some $\R$
can be thought of as solved problem if we have at our disposal the
\emph{constructive} representation for \emph{both} of these
functions. The `constructive' means here that they can be
manipulated analytically (differentiation, integration, etc) like
rational or elliptic functions.

First examples of algebraic solutions to equation \eqref{p6} were
obtained by R.~Fuchs in   work \cite{fuchs2} which is less known
than his 1907 work in Math.~Annalen with the same title as
\cite{fuchs2}. These solutions correspond to Picard's case of
parameters. Since the late 1990's list of algebraic solutions to
\Psix\ came into rapid growth thanks to works by Dubrovin, Mazzocco,
Kitaev, Boalch, Vid\a={u}nas, Hitchin himself, and others. Complete
reference list to this topic would be rather lengthy (see, \eg,
\cite{mazz, kitaev,boalch}) but recent work \cite{lisovyy} is
already devoted to classification results concerning \emph{all} the
algebraic solutions; a relevant analogy here is the `nonlinear
Schwarz's list' \cite{boalch2}. The same works
\cite{lisovyy,boalch2} contain also exhaustive references.

In this work we show that  correlation between algebraic
Picard--Hitchin solutions of the $\wp$-form \eqref{p6wp} and, on the
other hand, uniformizing Fuchsian equations leads to the new and
infinite families of explicitly parameterized algebraic curves
together with that constitutes the base of analysis on Riemann
surfaces: Abelian integrals and their differential calculus.  It is
worthy of special emphasis that while examples of uniformizing
functions (Hauptmoduln, \ie, principal moduli in Klein's terminology
\cite{fricke}) are known, few as they are, the state of the art as
to the integrals is still that no uniformizing $\tau$-representation
for Abelian integrals is known for any $\R$ of genus $g>1$. However
the detailed constructions of integrals will be postponed to a
separate continuation of this work since they lead to many other
interesting consequences.

\subsubsection{$\tau$-representations to the Picard--Hitchin class.
New $\vartheta$-constants}

The most interesting feature of the proposed family is that it is
infinite, uniformly  describable and fundamental group
representation of corresponding orbifolds (Riemann surfaces)
coincides with the automorphism group of a function field generator,
that is Painlev\'e function $y(\tau)$ itself. Moreover, the
computational apparatus involves a \emph{new} kind of
theta-constants and we develop  the differential calculus not only
for meromorphic objects on our $\R$'s (functions and differentials
as theta-ratios) but for these everywhere holomorphic functions as
well.

Picard--Hitchin's class lies in the base of our constructions but
extends further. Most nice examples we exhibit in the text are the
hyperelliptic curves  since they have numerous applications.
Although explicitly describable curves appear abundantly, we do not
touch on the classification problems; rather we stress the
\emph{ways of getting the analytic results}; a key-note of the work.
As a consequence, many problems accompanying the uniformization
theory  have got an explicitly solvable form.

\subsubsection{Consequences and applications}

Eq.~\eqref{p6} is rich in consequences going beyond the `proper
\Psix-theory' \cite{iwasaki,conte,CM,gromak,okamoto} and one of the
remarkable facet is the fact that all the ensuing results come from
a `physical' equation---the \Psix-equation. This explains the
presence of a large number of examples in text. Moreover, \Psix\ is
a first example of the natural origin where not only do algebraic
dependencies appear in a parameterized form, but we also have an
infinite family of such examples. All the $\mathcal{P}$-equations
have the extensive (math)physical literature and, on the other hand,
general theory of uniformization has long experienced the lack of
nontrivial illustrative applications. Even computational problems in
the theory are  not trivial and these received recently increased
attention in the context of another (aside from integrals and
functions) fundamental object: the Schottky--Klein prime function;
see \cite{crowdy1,crowdy3} and references therein. Although these
works are concerned with multiconnected versions of $\R$'s, the main
difficulties focus on computational and analytic aspects
anyway\footnote{D.~Crowdy notices properly in recent survey
\cite[p.~T206]{crowdy2}: `While many abstract theoretical results
exist, the mathematical literature has been lacking in
\emph{constructive} techniques for solving such problems~\..., the
mathematics applicable to such problems has widespread relevance to
all sorts of other problems in mathematical physics, classical
physics, integrable systems theory and the theory of ordinary and
\textsc{pde}s'.}.

As for applications, it is worth noting that inversions of Fuchsian
equations we consider give us the explicit and global
representations for canonical maps between moduli spaces of
one-parameter families of Calabi--Yau mirror pairs $(X, X^*)$
(IIA/IIB string duality), which are better known as the famous
mirror maps, realized through the Picard--Fuchs equations. See, \eg,
\cite{yau,yau2} for references source and \cite{hoeij} for a
computational application to the Ising model. On the other hand,
universal orbifold Hauptmoduln (see Sect.~\ref{universal} for
explanation what `universal' means here) parameterize
\emph{arbitrary $\R$'s of nonzero genera} and conversely, any
Riemann surface (punctured or compact) is necessarily described by
the zero genus Hauptmodul(n) \cite{br}. In particular, this entails
that there is actually no need to take, as a starting point, the
hyperbolic Poincar\'e polygons of genera $g>1$ because they
correspond to Fuchsian equations with algebraic coefficients. Among
other things, monodromies of such equations have, apart from
non-commutativity, very complicated geometric structure and their
hyperbolic generators may not be recognized from a local analysis.

\subsection{Outline of the work}

The article is organized as listed in Contents and consists
structurally of three parts; except for background
section~\ref{background}. The first part
(sects.~\ref{revis}--\ref{groups}) deals with uniformization of
Picard--Hitchin's algebraic solutions. The second part
(sects.~\ref{towers}--\ref{Fuchs}) is devoted to consequences,
relationships to other important equations, and examples. The third
part (sects.~\ref{con}--\ref{final}) sums up the two previous ones
from the general differential viewpoint with further emphasis on
toroidal covers and Abelian integrals.

In Sect.~\ref{background} we fix notation and briefly expound some
background information on uniformization theory: Fuchsian equations,
their monodromies, inversion problems, and transformations between
equations. We also introduce a convenient object---the meromorphic
$\bo{\mathfrak D}$-derivative---and the notion of universal
uniformization. At the heart of constructing the new uniformizable
curves is a simple lemma on transformations (Lemma~\ref{L1}).

In Sect.~\ref{revis} we recall the Painlev\'e substitution and
inspect some little-known facts about Picard's solution.

In Sect.~\ref{hit} we simplify the original parametric form of
Hitchin's solution and present results in the language of
uniformizing theory (Sect.~\ref{unihit}): nontrivial Fuchsian
equations and uniform representation for functions in terms of new
kind $\theta$-constants. These objects, along with the uniformizing
functions and Abelian differentials, admit the effective and
\emph{closed} differential computations. These technicalities are
expounded in Sect.~\ref{diffcalc}.

Section~\ref{groups} is devoted to  transformation groups. We first
tabulate the general transformation rules for $\theta$-functions  in
Sect.~\ref{groupstheta} and then, on the basis of this, in
Sect.~\ref{groupsPicard} we describe uniformizing groups for
Picard's curves.

Section  \ref{towers} consists basically of examples. One exhibits
ways of generation of uniformizable curves being no solutions to the
\Psix-equation; we call this a \emph{tower of curves}. Among these
are some non-hyperelliptic curves and three classical hyperelliptic
examples (sects.~\ref{Hyper} and \ref{shwarz}). We conclude this
section with the general recipe of getting the formulae.

It is notable that the Picard--Hitchin class has links with non
Picard--Hitchin's algebraic solutions to \Psix. The latter curves
also produce solvable Fuchsian equations. This is the subject matter
of Sect.~\ref{nonPH}.

In Sect.~\ref{Fuchs} we show that Picard--Hitchin's solutions have
direct links with many important Fuchsian equations. For example, a
simple Picard's solution yields an equation from the famous
Chudnovsky list of the four Heun equations (Sect.~\ref{Heun}). We
give some explanations as to how these equations can be related to
each other. The same Picard's solution is related, through
corresponding Fuchsian equation, to the famous Ap\'ery linear \ode;
we discuss this fact at greater length in Sect.~\ref{apery}. All
these facts allow us to write down the explicit
$\tau$-representations for associated inversion problems.  Although
the subsequent section (Sect.~\ref{further}) is devoted to further
illustration of the `$\theta$-apparatus' (including the nice
hyperelliptic example $z^2=x^6-1$) and relationships between
Hauptmoduln and solutions of Fuchsian equations, the main purpose of
this section is to anticipate an important generalization which will
be expounded in  the next section.

Section \ref{con} is devoted to differential structures on $\R$'s
and, in particular, to systematization of results of
sects.~\ref{diffcalc} and \ref{further}. We explain that, besides
the functions and differentials, one should introduce the notion of
analytic connection on $\R$. Sections~\ref{diffprop} and \ref{anal}
show how this object is described by means of certain \odes\ through
the unique fundamental scalar (automorphic function). We present the
regular recipe of getting such \odes\ and exhibit nontrivial
examples and exercises.

In final section  (Sect.~\ref{final})  we sketch a relationship
between the preceding material and transcendental (solvable)
Fuchsian equations on tori. Construction of such equations results
from the fact that the majority of our algebraic curves cover
elliptic tori. Since the toric coordinate has an (Abelian) integral
nature the integrals themselves should also be involved into the
theory.

\section{Background material\label{background}}

\subsection{Uniformization, Schwarzians, and Fuchsian
equations\label{nutshell}}

The classical language to describe uniformization of Riemann
surfaces of finite genera is the linear differential equations of
Fuchsian class \cite{ford}. If $x=\chi(\tau)$ is a generator of the
function field of meromorphic automorphic functions on a Riemann
surface $\R$ of some algebraic curve \eqref{F} then the global
uniformizing parameter $\tau$ is determined as a quotient
\begin{equation}\label{ratio}
\tau=\frac{\Psi_1(x)}{\Psi_2(x)}
\end{equation}
of two linearly independent solutions to the certain Fuchsian
equation of 2nd order \cite{ford,lehner}
\begin{equation}\label{fuchs}
\Psi_{\!\!\mathit{xx}}=\frac12\,\bcQ(x,y)\,\Psi\,.
\end{equation}
Equation of the same form
$\psi_{\mathit{yy}}=\frac12\,\widetilde\bcQ(x,y)\,\psi$, where
$\psi=\sqrt{\smash[bt]{y_x}}\:\Psi$, determines the $\tau$ through
the second function $y=\xi(\tau)$. That equation \eqref{fuchs} is of
Fuchsian class implies that function $\bcQ(x,y)$ is bound to be
rational. This function (or $\widetilde\bcQ$) completely determines
all the analysis on $\R$. The ratio $\tau$ itself, as a function of
$x$, is a solution of nonlinear non-autonomous \ode\ of 3rd order
\cite[p.~22]{ford}
\begin{equation}\label{multi}
\frac{\tau_{\!\mathit{xxx}}}{\tau_x}-\frac32\,
\frac{\tau_\mathit{xx}^2}{\tau_x^2}=-\bcQ(x,y)\,,
\end{equation}
better known as the  Schwarz equation, and left hand side of this
equation is traditionally designated by the Schwarz derivative
symbol \cite{schwarz}: $\{\tau,x\}=-\bcQ(x,y)$.

Solutions of equations \eqref{fuchs} and \eqref{multi} are
essentially multi-valued functions of variable $x$. This
multi-valuedness is described by a transformation group and the
group itself is nothing but matrix $(2{\times}2)$-representation of
the monodromy group $\G_x$ of equation \eqref{fuchs} \cite{ford}:
\begin{equation}\label{ma}
\G_x\colon\quad
\mbig[7](\,\begin{matrix} \Psi_1\\\Psi_2\end{matrix}\mbig[7])
\mapsto
\mbig[7](\begin{matrix} a&b\\c&d\end{matrix}\mbig[7])
\mbig[7](\,\begin{matrix} \Psi_1\\\Psi_2\end{matrix}\mbig[7])\,.
\end{equation}
This transformation entails  the main property of function
$x=\chi(\tau)$, namely, the property of being automorphic:
$$
\chi\Big(\Mfrac{a\,\tau+b}{c\,\tau+d}\Big)=\chi(\tau)\quad
\hence\quad \text{Aut}\;\chi(\tau)\FED\G_x\,.
$$

Let this $\chi(\tau)$, as analytic function of $\tau$, be globally
single-valued in the domain of its existence\footnote{The domain
$\mathbb{D}$ and property of being single-valued crucially depend on
parameters of function $\bcQ$. Those parameters that do not affect
the local analysis but determine the global domain and the global
single-valuedness (therefore topology of $\G_x$ and genus) are
called accessory parameters \cite{fricke,ford,yoshida,bers}. We
shall call these parameters \emph{correct} if $x=\chi(\tau)$ is a
single-valued function everywhere in $\mathbb{D}$.} $\mathbb{D}$ and
Poincar\'e polygon \cite{fricke,lehner} for monodromy $\G_x$ has
finite topological genus. We then may think of $\chi(\tau)$ as a
finite order meromorphic function on factor $\mathbb{D}/\G_x$. The
latter is generally an orbifold \cite{yoshida} and, upon
compactification (if required), may be turned to a Riemann surface
of some algebraic curve \eqref{F}, possibly sphere
$\mathbb{P}^1(\mathbb{C})$. It is known that this construction can
always be realized by a Kleinian group with an additional condition
that the group has an invariant circle and determines thereby a
Fuchsian group of 1st kind \cite{ford,lehner}. We normalize this
circle to be the real axis $\mathbb{R}$ and universal cover (where
$\tau$ `lives') to be the upper-half plane $\Hp\ni\tau$, that is
$\Im(\tau)> 0$;  the domain $\mathbb{D}$ thus becomes $\Hp$.

If Poincar\'e domain for the single-valued $\chi(\tau)$  has some
punctures then the local behavior of $\chi(\tau)$ near to them is
well-defined, \ie, has a local parameter without ambiguities. Hence,
we may do a (unique) compactification $\overline{\HG}$ of this
polygon \cite{bers} and (analytic) functions on that compactified
object become equally well as functions on pure hyperbolic Riemann
surfaces without punctures. See Remarks~\ref{R1} and \ref{R5}
further below for additional explanations and Sect.~\ref{Hyper} for
examples.

With the exception of accessory parameters problem  the most
important problems in the field are 1) explicit solutions to
\eqref{fuchs} and 2) explicit representation for inversions of the
ratio \eqref{ratio}. All the currently known solutions to the first
of these problems are reduced to the hypergeometric functions and
triangle groups \cite{fricke}; though there are  curves (Shimura
curves) with their Fuchsian equations \eqref{fuchs} having no such a
type of reduction (see, \eg, \cite{chud2,krammer,elkies} for
explicit formulae).

As for the second problem (the inversion problem) the number of
solvable examples falls far short of the first one and is limited
only by particular triangle groups, namely, groups commensurable
with the full modular group
$\bo\Gamma(1)\DEF\mathrm{PSL}_2(\mathbb{Z})$. Curiously, even the
most famous curve $x^3y+y^3z+z^3x=0$ was uniformized by Klein
\cite{levy} not through its `native' and famous hyperbolic
(2,3,7)-triangle but through the theta-constants associated with
group $\bo\Gamma(7)$, \ie, matrices congruent to the unity modulo~7
\cite{levy}. This point is a manifestation of the fact that no one
representation for uniformizing function associated with any
non-modular Fuchsian equation is known hitherto. In this context, to
the best of our knowledge, mentioning the Abelian integrals is
almost absent \cite[p.~551]{dalz1}, \cite{dalz2} and some compatible
constructions `integrals + modular groups' \cite{knopp1,knopp2} are
still far from being a theory.

\subsection{Meromorphic derivative}

Algebraic functions like $x$, $y$, etc are not only meromorphic
single-valued objects on $\R$'s. Complete analysis should
necessarily include meromorphic Abelian differentials and their
integrals as well; see also Sect.~\ref{diffprop}. Rather than
manipulate with non-autonomous equations \eqref{multi} and
multi-valued inversions of multi-valued objects it is convenient to
invert Schwarz's derivative $\{\tau,x\}$ to the object
$\{x,\tau\}/\dot x^2$ and to handle the \emph{meromorphic
derivative} $\bo{\mathfrak D}$ \cite{br3}:
$$
\bo{\mathfrak D}\colon\qquad [x,\tau]\DEF
\frac{\dddot{\smash[b]{x}}}{\dot x^3}-\frac32\,\frac{\ddot x^2}
{\dot x^4}\,,
$$
where dots above the symbols, as always in the sequel, stand for
derivatives with respect to $\tau$. For the reasons above,  there is
no point in distinguishing Fuchsian linear or Schwarz's nonlinear
equations
\begin{equation}\label{MD}
[x,\tau]=\bcQ(x,y)\,.
\end{equation}
Hence, in order to avoid lengthening terminology we shall use  the
standard notions---regular singularities, branch places, Fuchsian
exponents, et cetera \cite{ford,yoshida}---for both kinds of
equations. In particular, we shall use as synonyms the notions
automorphisms $\text{Aut}$, monodromies, and groups like $\G_x$ and
refer to linear equations \eqref{fuchs} and their nonlinear Schwarz
varieties \eqref{MD} merely as Fuchsian equations. The famous
theorem of Klein, Poincar\'e \& Koebe \cite{fricke,bers} guaranties
the availability of what is called the unique Fuchsian monodromy of
1st kind for equation \eqref{MD}, that is uniqueness of function
$\bcQ(x,y)$ with a given set of $x$-singularities $\{E_k\}$.

Throughout the paper we consider as equivalents Fuchsian equations
in the normal Klein form \eqref{fuchs} and equations with the
structure
\begin{equation}\label{pqPsi}
\psi_{\mathit{xx}}^{}+p(x,y)\,\psi_x^{}+q(x,y)\,\psi=0
\end{equation}
because transition between these forms is achieved by the known
linear transformation
$$
\Psi=\psi\cdot\exp\frac12\Smaller{\ds\int}\! p\:dx\,.
$$
Clearly, it has no effect on the ratio \eqref{ratio}.

\subsection{Universal uniformization and transformations between
Fuchsian \odes\label{universal}}

If the structure of $\bcQ$ is such that in the neighborhood of some
point $x=E$ we have
$$
[x,\tau]=-\frac12\frac{1}{(x-E)^2}+\cdots\,,
$$
then the local behavior of function $\chi(\tau)$ will be exponential
\cite[Ch.~5]{yoshida}, \cite{ford}. This being so, arbitrary
algebraic ramification $y(x)\sim (x-E)^q+\cdots$ ($q\in\mathbb{Q}$)
is transformed into the locally single-valued dependence $y(\tau)$.
This point motivates a definition for the universal uniformization
through punctures.

\medskip
\noindent
{\bf Definition.} \emph{The meromorphic automorphic
function $x=\chi(\tau)$ on $\Hp$ is said to be the universal
uniformizing function for the set of points $\{E_k\}$ if it has the
exponential behavior
$$
x=E+a\exp\left(\mfrac{-\pi\ri}{\tau-\tau_0^{}}\right)+
\cdots
$$
in neighborhoods of $E$'s as $\tau\to \tau_0^{}+0\,\ri$ under
$\tau_0^{}\in\mathbb{R}$. If $\tau_0^{}\to\ri\infty$ we define}
$$
x=E+a\exp(\pi\ri\,\tau)+\cdots\,.
$$

In other words, local monodromy $\G_x$ for universal uniformizing
function is determined by the parabolic singularities $E_k$ in
\eqref{fuchs} and implies the exponents above. The definition does
not forbid to have a non-parabolic, \ie, conical, singularity. It
follows that any algebraic function of $x$, say \eqref{F}, with
\emph{arbitrary} ramifications only at points $x=E_k$ becomes a
single-valued function of $\tau$. Correspondence between the
exponential behavior, that is puncture, and the $\bo{\mathfrak
D}$-object  is perhaps most easily clarified by observing that
meromorphic derivative of a function can be represented through the
one of its logarithm:
$$
[x,\tau]=\frac{1}{(x-e)^2}\!\left(\big[\ln (x-e),\tau\big]-
\frac12\right)
$$
with arbitrary $e$; logarithm is `swallowed up' by the exponent
above. The classical and simplest example is a uniformization of a
3-punctured sphere:
\begin{equation}\label{k2}
[x,\tau]=-\frac12\,\frac{x^2-x+1}{x^2(x-1)^2}\,.
\end{equation}
Solution of this equation is the $x$-function in Painlev\'e
substitution; see formula \eqref{subswp} further below.

\begin{remark}\label{R1}
An interrelation needs to be understood between parametrizations,
uniformization, and the curve itself. Analytically,  curve is not an
invariant object since we may do  birational transformations. It
therefore has no punctures or conical singularities. The only
invariant object is a Riemann surface $\R$ with a pure hyperbolic
system of generators representing its fundamental group
$\pi_{\!1}(\R)$. On the other hand, uniqueness of compactification
entails that complex analytic objects---generators $x(\tau)$,
$y(\tau)$ of a function field and Poincar\'e domain
$\mathbb{D}$---have been subordinated to the only fundamental
analytic property: the global single-valuedness with due regard for
corresponding factor topology. As for the differential apparatus, it
is local at all and does not `sense' the global and conformal
requirements. It and many other ingredients of the theory do need
just an analytic equivalence; not the global conformance. These are
not one and the same because analytic equivalence of $\R$'s (curves)
is weaker than the conformal one. In this respect uniformization
with punctures is not inferior than pure hyperbolic
one\footnote{Indeed, considering the object/sphere $x^2+y^2=4$ we do
not think of its two parametrizations
$\{x=2\sin\tau,\,\,y=2\cos\tau\}$ and $\{x=\tau+\tau^{\sm1},\ri
y=\tau-\tau^{\sm1}\}$ as one is `better' and other is `worse'. It is
clear that the obvious transformation $\tau\mapsto \re^{\ri\tau}$
realizes a translation between  these `punctured' and
`non-punctured' uniformizations. The remarkable fact is that there
exists an analog of that transition for higher genera and it is
described by \odes; the first explicit instance for the case $g=2$
has been exhibited in work \cite{br3}.}. Moreover, punctures are not
necessary conditions to construct universal uniformization since
there exist `non-punctured' universal ones. Every compact Riemann
surface $\R$ does indeed have a naturally associated orbifold
$\bo{\mathfrak{T}}$ \emph{without} punctures and the very pure
hyperbolic representation of $\pi_{\!1}(\R)$ turns out to be merely
a subgroup of $\pi_{\!1}(\bo{\mathfrak{T}})=\G_x$, where
$\G_x=\text{Aut}(x)$, and function $x=\chi(\tau)$ is automorphic
with respect to representation of this $\pi_{\!1}(\R)$. Topological
arguments show that group $\G$ uniformizing the curve \eqref{F} and
representing $\pi_{\!1}(\R)$ is to be taken as $\G=\G_x \cap\G_y$.
This
is the subject matter %(including extended bibliography)
of work \cite{br}.
%All this,
%including extended bibliography,
%and some further explanations of universal
%uniformization is expounded at greater length in work \cite{br}.
It is of interest to remark here that Fuchsian equations for
subgroups/curves have not got % (strangely enough)
to the second  volume of the monumental Fricke--Klein treatise on
automorphic functions \cite{fricke} (devoted to a function theory),
whereas Picard--Hitchin's and many other modular equations provide
examples of such constructions.
\end{remark}

Meromorphic derivative $\bo{\mathfrak D}$ can be calculated for any
object on $\R$. If this object is an element of function field then
it is algebraically related to any other one. This algebraic
relation, considered as a change of variable, transforms Fuchsian
equations one into the other.

\begin{lemma}\label{L1}
If \,$z=R(x)$ is any function of $x$ then equation
$[x,\tau]=\bcQ(x,y)$
implies that
\begin{equation}\label{*}
[z,\tau]=[R(x),x]+\frac{1}{R_x^2}\,\bcQ(x,y)\,.
\end{equation}
\end{lemma}

This is of course the $\bo{\mathfrak D}$-version of the known
transformation law for Schwarz's derivative of a function
composition $\tau\circ \mu$. If $\tau=f(\mu)$ and $\mu=g(z)$ then we
have
\begin{equation}\label{tensor}
\{\tau,\mu\}\,d\mu^2+\{\mu,z\}\,dz^2=\{\tau,z\}\,dz^2\,.
\end{equation}

In spite of seemingly triviality, this lemma has a fundamental
meaning because algebraic curves may form \emph{towers} and
integrability of their Fuchsian equations is in effect an
integrability of a \emph{single} equation (see Sect.~\ref{towers}).
Global parameters $\tau$'s for all of these equations/curves are the
one common $\tau$ and we search for relations between these curves.
In our case, these relations are the algebraic/rational ones to
Hauptmodul $x=\chi(\tau)$; the problem consists in finding these
substitutions and representations for various (Haupt) Moduln.

In the following, we shall exhibit such results coming from the
Picard--Hitchin solutions to \Psix. For example, nontrivial
substitution may turn the equation \eqref{*} with complicated
function $\bcQ(x,y)$ into a simple equation of the form
$[z,\tau]=\bcQ(z)$ even though the curve \eqref{F} has a nontrivial
genus (see Sect.~\ref{7.3}). Of course, our parametrizations
correspond to finite covers of the punctured spheres and orbifolds
since the group $\bo\Gamma(1)$ and its subgroups like $\bo\Gamma(2)$
form presently the only class for which explicit inversions of
\eqref{ratio} are known. However, as pointed out above, we get a
massive extension of the classical family of Jacobi's
$\vartheta$-constants and Dedekind's eta-function.

\subsection{Notation\label{not}}

Picard--Hitchin's solutions involve nontrivial combinations of
Jacobi's and Weierstrass's functions and we use intensively many of
their properties without explicit mentioning. Among enormous
literature on this subject, in most of cases Schwarz's collection of
Weierstrass's and Jacobi's classical results \cite{we2} is by no
means lacking and the four volume set by Tannery \& Molk
\cite{tannery}, as a formulae source, hitherto contains most
exhaustive information along these lines.

We use the four Jacobi's functions $\theta_k$, introduced by Hermite
as $\theta$-functions with characteristics \cite[p.~482]{hermite},
in the following definition \cite{weber}:

\begin{equation}\label{theta}
\theta\AB{\alpha}{\beta}(z|\tau)=
\sideset{}{_k}\sum\limits_{\sm\infty}^{\infty}\!
\re_{\strut}^{\pi\ri \s[2](k+
\frac\alpha2\s[2])^2\tau+
2\,\pi\ri\s[2](k+\frac\alpha2\s[2])
\s[2](z+\frac{\smash{\beta}}{2}\s[2])}\,.
\end{equation}
Therefore, $\theta_1=-\theta\AB{1}{1}$, $\theta_2=\theta\AB{1}{0}$,
$\theta_3=\theta\AB{0}{0}$, $\theta_4=\theta\AB{0}{1}$. Nullwerthe
of $\theta$'s, termed usually the $\vartheta$-constants, are the
values of $\theta(z|\tau)$ under $z=0$, that is $\vartheta_k\DEF
\vartheta_k(\tau)=\theta_k(0|\tau)$. We introduce the fifth and
independent object $\Dtheta$ as a derivative of the
$\theta_1$-series:
$$
\Dtheta(z|\tau)= \pi\,\re^{\frac14\pi\ri\tau}_{\mathstrut}
\sideset{}{_k}\sum\limits_{\sm\infty}^{\infty}\! (-1)^k\,
(2k+1)\,\re^{(k^2+k)\,\pi\ri\tau}\, \re^{(2k+1)\,\pi \ri z}\,.
$$

The standard Weierstrassian functions
$\sigma,\zeta,\wp,\wp'(z|\omega,\omega')$ correspond to the set of
half-periods $(\omega,\omega')$ \cite{we2,weber,tannery}. By virtue
of homogeneous relations, say
$\alpha^2\,\wp(\alpha\,z|\alpha\,\omega,\alpha\omega')=
\wp(z|\omega,\omega')$, we may always put any half-period to unity
and, as throughout the paper, handle functions like
$\wp(z|1,\tau)\FED\wp(z|\tau)$, etc, where $\tau$ stands for the
ratio $\omega'\!/\omega$. If $\Im(\tau)<0$, then one swaps around
$\omega$ and $\omega'$. Period of the canonical meromorphic elliptic
integral $\zeta$ is usually denoted as
$\eta(\tau)\DEF\zeta(1|\tau)$. Translation of Weierstrassian
functions into the $\theta$-language is realized by means of the
following formulae \cite{we2,tannery}:
\begin{equation}\label{wp}
\begin{array}{c}
\ds\wp(2z|\tau)=\frac{\pi^2}{12}
\bigg\{\vartheta_3^4(\tau)+\vartheta_4^4(\tau)+3\,\vartheta_3^2(\tau)
\,\vartheta_4^2(\tau)\,
\frac{\theta_2^2(z|\tau)}{\theta_1^2(z|\tau)}
\bigg\},\\[3ex]
\ds\zeta(2z|\tau)=2\,\eta(\tau)\,z+\frac12\,\frac{\Dtheta(z|\tau)}
{\theta_1(z|\tau)}\,,\qquad
\wp'(2z|\tau)=-\pi^3\,\ded^9(\tau)\,\frac{\theta_1(2z|\tau)}
{\theta_1^4(z|\tau)}\,,
\end{array}
\end{equation}
where $\ded(\tau)$ is the function of Dedekind:
$$
\oper\eta\ded(\tau)= \re^{\frac{\pi\ri}{12}\,\tau}_{\mathstrut}\,
\sideset{}{_k}\prod\limits_1^\infty
\big(1-\re^{2\,k\,\pi\ri\tau}\big)=
\re^{\frac{\pi\ri}{12}\,\tau}_{\mathstrut}
{\sideset{}{_k}\sum\limits_{\sm\infty}^{\infty}}
\!(-1)^k\,\re^{(3k^2+k)\,\pi\ri\tau}\,.
$$
It is differentially and algebraically related to the
$\eta,\vartheta$-constants:
\begin{equation}\label{difalg}
\frac{1}{\ded}\,\frac{d\ded}{d\tau}=\frac{\ri}{\pi}\,\eta\,,\qquad
2\,\ded^3(\tau)=\vartheta_2(\tau)\,\vartheta_3(\tau)\,
\vartheta_4(\tau)\,.
\end{equation}

\section{Picard's solution, revisited\label{revis}}

Contemporary mentions of this solution refer usually  to \eqref{p6}
however \ode\ derived by Picard himself differs from Painlev\'e form
\eqref{p6}.   As said above, Picard was not concerned with equation
\eqref{p6} or some of its particular case. At the end of his
m\'emoire \cite{picard} he offered as an example  the une \'equation
diff\'erentielle curieuse with fixed critical points satisfied by
Jacobi's elliptic sinus $\mathrm{sn}(a\,\omega+b\,\omega';k)$
considered as function of Legendre's modulus $k^2=x$. Picard denoted
this function by $u(x)$ and deduced the equation
\begin{equation}\label{u}
\begin{aligned}
\frac{d^2u}{dx^2}&-\Big(\frac{du}{dx}
\Big)^{\!2}\,\frac{u\,(2\,x\,u^2-1-x)}{(1-u^2)(1-x\,u^2)}\\
&+\frac{du}{dx}\bigg[ \frac{u^2-1}{(1-x)(1-x\,u^2)}+\frac1x\bigg]-
\frac14\frac{u\,(1-u^2)}{x\,(1-x)(1-x\,u^2)}=0
\end{aligned}
\end{equation}
(in original  paper \cite{picard}  on p.~298 the multiplier
$\frac14$ in front of last term was missing).

Transformation between Picard's and Painlev\'e equations
\eqref{p6wp}, \eqref{u} can be derived from the substitution of
Painlev\'e $(x,y)\mapsto (z,\tau)$ turning equation \eqref{p6} into
\eqref{p6wp}. Indeed, converting Painlev\'e formulae
\cite[p.~1117]{painleve} into the theta-functions, we get the
substitution \cite{babich}
\begin{equation}\label{substheta}
x=\frac{\vartheta_4^4(\tau)}{\vartheta_3^4(\tau)}\,,\qquad
y=-\frac{\vartheta_4^2(\tau)}{\vartheta_3^2(\tau)}\,
\frac{\theta_2^2\big(\frac12 z\big|\tau
\big)}{\theta_1^2\big(\frac12z\big|\tau \big)}
\end{equation}
and its  $(\vartheta,\wp)$-equivalent \cite{manin}
\begin{equation}\label{subswp}
x=\frac{\vartheta_4^4(\tau)}{\vartheta_3^4(\tau)}\,,\qquad
y=\frac13+\frac13\frac{\vartheta_4^4(\tau)}{\vartheta_3^4(\tau)}
-\frac{4}{\pi^2} \frac{\wp(z|\tau)}{\vartheta_3^4(\tau)}\,.
\end{equation}
Hence relation between Picard's $u$ and Painlev\'e' $y$ is nothing
but the known relation between functions $\wp$ and $\mathrm{sn}$:
$$
\wp(z+\tau|\tau)-\e'(\tau)=\frac{\pi^2}{4}\,\vartheta_2^4(\tau)\cdot
\mathrm{sn}^2\Big(\mfrac{\pi}{2}\vartheta_3^2(\tau)\,z;k \Big)
$$
and therefore
\begin{equation}\label{change}
u\Big(\Mfrac1x\Big)= \sqrt{y(x)}\,.
\end{equation}
The first equation in \eqref{subswp} is invertible and inversion
itself involves  Legendre's functions $\ellK$ and $\ellK'$
\cite{tannery}:
\begin{equation}\label{tau}
\tau=\ri \frac{\ellK(\!\sqrt{x})}{{\ds\ellK'}(\!\sqrt{x})}\,.
\end{equation}
The quantities $\ellK$ and $\ellK'$ can be considered as
hypergeometric series
\begin{equation}\label{kk'}
\ellK(\sqrt{x})=\frac{\pi}{2}\cdot
{}_2F_1\Big(\Mfrac12,\Mfrac12;1\Big|x \Big)\,,\qquad
\ellK'(\sqrt{x})=\frac{\pi}{2}\cdot
{}_2F_1\Big(\Mfrac12,\Mfrac12;1\Big|1-x \Big)
\end{equation}
or, depending on preference, as complete elliptic integrals
\cite{abramowitz}
$$
\ellK(\sqrt{x})=\int\limits_0^{\,\,1}\!\!\frac{d\lambda}
{\sqrt{(1-\lambda^2)(1-x\,\lambda^2)}}\,,\qquad
\ellK'(\sqrt{x})=\int\limits_{\sqrt{x}}^{\,\,1}\!\!\frac{d\lambda}
{\sqrt{(\lambda^2-x)(1-\lambda^2)}}\,.
$$
The relation \eqref{change}  leads to the following property of
Picard's case.

\begin{proposition}\label{P1}
Picard's solution of equation \eqref{p6} is a perfect square\,$:$
\begin{equation}\label{algPic}
y_{\sss\textsc{Pic}}^{\mathstrut}=-\sqrt{x}\;
\frac{\theta_2^2\mbig(
A\frac{\ellK(\!\sqrt{x})}{\ellK'(\!\sqrt{x})}+B\mbig|\ri\frac{
\ellK(\!\sqrt{x})}{\ellK'(\!\sqrt{x})}\mbig)}{\theta_1^2
\mbig( A\frac{\ellK(\!\sqrt{x})}{\ellK'(\!\sqrt{x})}+
B\mbig|\ri\frac{\ellK(\!\sqrt{x})}{\ellK'(\!\sqrt{x})}\mbig)}\,,
\end{equation}
where $A$ and $B$ are free constants.
\end{proposition}

This is a special feature of Picard's parameters since square root
of arbitrary solution to equation \eqref{p6} contains no movable
ramifications and satisfies another equation with the Painlev\'e
property, that is equation \eqref{u}.  In this regard Hitchin's
solution is much more nontrivial (see Sect.~\ref{hit}). The
$\wp$-form of Picard's solution was considered comprehensively by
Fuchs \cite{fuchs2} and reinspected by Mazzocco in \cite{mazz}.

On the other hand, we know that Gambier's list of fifty
transcendents \cite{gambier} is complete and, contrary to
\eqref{change},  there is bound to be a change of independent
variable $x\mapsto z=f(x)$ transforming \eqref{u} into the one of
equations from this list. It is not difficult to see that it can be
only Picard's case of equation \eqref{p6}. See the number
(\oldstylenums{49}) on p.~17 in \cite{gambier} and, concerning the
change itself, the case III on p.~323 in \cite{ince}; the latter
transforms  into the case viii on p.~326 in the same place. Carrying
out calculations and  gathering intermediate substitutions we get
the following symmetry.

\begin{proposition}\label{P2}
Let $\alpha=\beta=\gamma=\delta=0$. Then the transformation
$$
y(x)\mapsto
\left\{
\frac{y\mbig[7](\Big(\mfrac{\sqrt{x}-1}{\sqrt{x}+1}\Big)^2\mbig[7])+
\mfrac{\sqrt{x}-1}{\sqrt{x}+1}}
{y\mbig[7](\Big(\mfrac{\sqrt{x}-1}{\sqrt{x}+1}\Big)^2\mbig[7])-
\mfrac{\sqrt{x}-1}{\sqrt{x}+1}}\right\}^{\!2}
$$
is a function automorphism of equation \eqref{p6} preserving the
perfect square.
\end{proposition}

Transformations of such a kind are the subject of a more general
theory of quadratic transformations and additional details on this
theory can be found in work \cite{kitaev}.

\section{Uniformization of Picard--Hitchin's curves\label{hit}}

\subsection{Hitchin's solution}%, revisited}

General group of invariance of the Picard--Painlev\'e property is a
birational group \cite{conte}. This transformation group includes
derivatives (see Okamoto's transformations \eqref{ok} further below)
and leads to Hitchin's solution which has  hitherto remained  most
nontrivial and `rich' solution of all those currently known. Its
parametric form is as follows \cite[Theorem~6]{hitchin}:
\begin{equation}\label{hitchin1}
\wp(z|\tau)=\wp(A\tau+B|\tau)+\frac12\, \frac{\wp'(A\tau+B|\tau)}
{\zeta(A\tau+B|\tau)-A\,\eta'(\tau)-B\,\eta(\tau)}
\end{equation}
with the same meaning of $A$, $B$, and $\tau$ as in
Proposition~\ref{P1}. In implicit form this solution was also
obtained by Okamoto \cite[p.~366]{okamoto}. In addition to
\eqref{hitchin1}, Hitchin suggested a theta-function form for his
solution which turned out to be rather complicated since it contains
the set of functions
$\vartheta,\vartheta',\vartheta''',\theta,\theta',\theta'',
\theta'''$.
Reproduce the solution in Hitchin's original notation on p.~33 of
\cite{hitchin}:
\begin{equation}\label{morazm}
\begin{aligned}
y(x)&=\frac{\vartheta_1'''(0)}{3\,\pi^2\vartheta_4^4(0)\,
\vartheta_1'(0)}
+\frac13\!\left(1+\frac{\vartheta_3^4(0)}{\vartheta_4^4(0)}
\right)\\
&\== +\frac{\vartheta_1'''(\nu)\,\vartheta_1^{}(\nu)-
2\,\vartheta_1''(\nu)\,
\vartheta_1'(\nu)
+4\,\pi\ri c_1^{}\,\big(\vartheta_1''(\nu)\,
\vartheta(\nu)-\vartheta_1'^2(\nu)
\big)}{2\,\pi^2\,\vartheta_4^4(0)\,\vartheta_1^{}(\nu)\,
\big(\vartheta_1'(\nu)+
2\,\pi\ri c_1^{}\,\vartheta_1^{}(\nu)\big)}\,,
\end{aligned}
\end{equation}
where $\nu=c_1^{}\tau+c_2^{}$. Below is simplification of Hitchin's
formulae resulting in a full parametrization of subsequent curves.

\begin{theorem}\label{T1}
General solution to Hitchin's case of equation \eqref{p6} has the
form
\begin{equation}\label{simp}
y=\frac{\sqrt{x}}{\theta_1^2}
\bigg\{\frac{\pi\,\vartheta_2^2\cdot\theta_2\,\theta_3\,\theta_4}
{\Dtheta+2\,\pi A\,\theta_1}-\theta_2^2\bigg\},
\end{equation}
where functions $\Dtheta$, $\theta_k$ are understood to be equal to
$\Dtheta,\,\theta_k\mbig[5](A\frac{\ellK(\sqrt{x})}
{\ellK'(\sqrt{x})}+
B\mbig[4]|\ri\frac{\ellK(\sqrt{x})}{\ellK'(\sqrt{x})}\mbig[5])$ and
$\vartheta_2\DEF\vartheta_2\mbig[5](\ri
\frac{\ellK(\sqrt{x})}{\ellK'(\sqrt{x})}\mbig[5])$.
\end{theorem}

\begin{pf}
Meromorphic functions on Riemann surfaces are expressed through
ratios of theta-functions and therefore the set of
$\theta$-functions in \eqref{morazm} is excessive. Based on
conversion formulae \eqref{wp} and duplication rules for
$\theta$-functions \cite{tannery} we obtain that \eqref{hitchin1}
$\dashrightarrow$ \eqref{simp}. First formula in \eqref{wp} explains
also the transition between \eqref{substheta} and \eqref{subswp}.
\hfill $\blacksquare$
\end{pf}

\begin{remark}\label{R2}
An important point here is the fact that function $\Dtheta$ should
play an independent part in the $\theta$-calculus along with  the
four functions $\theta_k$. An explanation of this role of $\Dtheta$
is as follows. Weierstrass' $\zeta$-function is the canonical
meromorphic integral on elliptic curves. It never reduces to
holomorphic or logarithmic integrals, or meromorphic elliptic
functions since it is related to them only through a derivative. On
the other hand $\zeta$-function is proportional to the ratio
$\Dtheta/\theta_1$.
\end{remark}

\subsection{Uniformization of Hitchin's  curves\label{unihit}}

Algebraic solutions \eqref{F} of Picard's class correspond to
constants \cite{picard,fuchs2,mazz}
\begin{equation}\label{nmN}
A\tau+B=\frac{\nu}{2N}\,\tau+\frac{\mu}{2N}\,,\qquad N\in
\mathbb{Z}\backslash\{0\}\,.
\end{equation}
By virtue of Theorem~\ref{T1} algebraic solutions to the Hitchin
case are effectively described as well as Picard's class.

\begin{theorem}\label{T2}
Let $\tau\in\Hp$, $\vartheta_k=\vartheta_k(\tau)$, and functions
$\Dtheta,\theta_k(A\tau+B|\tau)$ be taken with arguments
\eqref{nmN}. Then algebraic solutions \eqref{F} to the Hitchin case
of parameters in \eqref{p6} have the following parametrization by
single-valued functions\/$:$
\begin{equation}\label{param}
x=\frac{\vartheta_4^4}{\vartheta_3^4}\,,     \qquad
y=\frac{\vartheta_4^2}{\vartheta_3^2}\,\frac{\theta_2}{\theta_1^2}
\bigg\{
\frac{\pi\,\vartheta_2^2\cdot\theta_3\,\theta_4}{\Dtheta+\pi\ri
\frac{\nu}{N}\,\theta_1}-\theta_2 \bigg\}.
\end{equation}
\end{theorem}

\begin{pf}
That functions \eqref{param} satisfy some algebraic dependence
\eqref{F}  will be apparent from the fact that \eqref{nmN} generates
an algebraic family of Picard $P$ \eqref{algPic} and the latter is
related to Hitchin's solutions $H$ through the Okamoto
transformations \cite{okamoto,gromak}:
\begin{equation}\label{ok}
H=P+\frac{P\,(P-1)(P-x)}{x\,(x-1)\,P_x-P^2+P}\,,\qquad
P=H-\frac{H\,(H-1)(H-x)}{
x\,(x-1)\,H_x+\frac12\,H^2-x\,H+\frac12\,x}\,.
\end{equation}
These maps transform algebraic dependencies into the same ones.
\hfill $\blacksquare$
\end{pf}

\begin{corollary}\label{C1}
Let us fix $N$ in \eqref{nmN}.  Then every Picard's solution is a
rational function $P=R(x,H)$ on Hitchin's curve $F_N(x,H)=0$ and
conversely.  Picard's and Hitchin's solutions, along with any
algebraic solutions tied by a certain Okamoto transformation, define
isomorphic curves.
\end{corollary}

\begin{example}\label{E1}
Under $N=3$ and $(\nu,\mu)=(0,1)$ the birational isomorphism reads
as follows
$$
P=1-\frac{(x-1)^2}{(H-1)^2}\,,\qquad H=\frac32\,\frac xP-
\frac12\,P\,.
$$
\end{example}

In work \cite{br} we conjectured that theta-constants of the general
form $\theta\big(u(\tau)|\tau\big)$ can be relevant in
uniformization theory. As we have seen now this is so indeed and the
new family of theta-constants
$\theta\big(\frac{\nu}{N}\tau+\frac{\mu}{N}\big|\tau\big)$ comes
into play. The distinctive property of this family, as compared with
classical modular $\ded,\vartheta$-constants \cite{weber}, is an
availability of the special constant
$\Dtheta\big(\frac{\nu}{N}\tau+\frac{\mu}{N}\big|\tau\big)$ and its
role in differential closure. One more feature is the general
formula for all the family of uniformizing functions, which cannot
be said of many subfamilies of modular equations. For example,
algebraic dependencies between roots of Legendre's modulus
$k^2(\tau)$ have the form $F\big(\!\sqrt[\uproot{1}n]{k\,}(p\,\tau),
\sqrt[\uproot{1}m]{k\,}(s\,\tau)\big)=0$, which is easily written
down for any integer $n,m,p,s$, but uniformizing functions are known
explicitly only in cases $n,m=\{1,2,4\}$ (see the last paragraph in
Sect.~14.6.3 of \cite{bateman}). Among other things, the size of
coefficients in modular equations, as is well known, rapidly growths
as level increases \cite{weber}. We shall use the term
\emph{$\theta$-constant} (\emph{$\theta'$\!-constant}) if
$z$-arguments of the $\theta,\,\theta'$-functions \eqref{theta} are
some functions of the modulus~$\tau$.

\begin{theorem}\label{T3}
Let $y(\tau)$ be a uniformizing representation \eqref{subswp} of the
arbitrary algebraic Painlev\'e solution \eqref{F}. Then $y(\tau)$ is
single-valued and satisfies the Fuchsian \ode\/$:$
\begin{equation}\label{Qxy}
[y,\tau]= \big([F,x]-[F,y]\big)\,\pow Fy2 +3\,\frac{F_y}{F_x}
\![3]\left(\ln\!\frac{F_y}{F_x}\right)_{\![4]\mathit{xy}}\!
-\frac12\,\frac{x^2-x+1}{x^2(x-1)^2}\,
\frac{\pow{F}{y}{2}}{\pow{F}{x}{2}}\,,
\end{equation}
where $[F,x]$ and $[F,y]$ are computed as the  partial
$\bo{\mathfrak D}$-derivatives.
\end{theorem}

\begin{pf}
Algebraic dependence \eqref{F} and Lemma~\ref{L1} imply that
Fuchsian equations \eqref{MD} for automorphic functions $x(\tau)$
and $y(\tau)$ are not independent. The function $x(\tau)$ is an
elliptic modulus $k'^2(\tau)$; therefore it satisfies equation
\eqref{k2} and is a universal uniformizing function for the three
points $x=\{0,1,\infty\}$. Since Painlev\'e functions ramify only
over these points, $y(\tau)$ is single-valued. Applying
Lemma~\ref{L1} and substitution \eqref{F} to equation \eqref{k2}, we
get, after some simplification, the computational rule \eqref{Qxy}.
\hfill $\blacksquare$
\end{pf}

In what follows we shall show that equations \eqref{Qxy}, being
nontrivial Fuchsian ones with algebraic coefficients, are integrable
and have  computable monodromies.

In many respects, the simplicity of functions \eqref{param} leads to
that we have actually no need for equation of the curve itself. In
addition to Puiseux series, differentiations, plotting of graphs,
etc, this is especially true in applications wherein parametric
representation is the best suited form for implicit solutions.

\section{Differential calculus on Picard--Hitchin
curves\label{diffcalc}}

Complete analysis of meromorphic functions on Riemann surfaces must
contain Abelian differentials but Theorems~\ref{T2} and \ref{T3} do
not touch on these objects. Picard--Hitchin's curves (abbreviated
further to PH-curves) contain $\vartheta$-constants and some of
their differential properties (not all) are known \cite{weber}. The
following result has a technical characterization but is needed for
completeness of differential computations with theta-functions.

\begin{lemma}[\cite{br}]\label{L2}
Jacobi's $\vartheta$-constants form a differential ring over
$\mathbb{C}(\pi)$ upon adjoining the period of meromorphic elliptic
integral, that is Weierstrass' $\eta{:}$
\begin{equation}\label{var}
\begin{aligned}
\frac{d\vartheta_2}{d\tau}&=
\frac{\ri}{\pi}\bigg\{\eta+\frac{\pi^2}{12}\,
\big(\vartheta_3^4+\vartheta_4^4
\big)\bigg\}\vartheta_2\,,\\
\frac{d\vartheta_3}{d\tau}&=
\frac{\ri}{\pi}\bigg\{\eta+\frac{\pi^2}{12}\,
\big(\vartheta_2^4-\vartheta_4^4 \big)\bigg\}\vartheta_3\,,
\end{aligned}\qquad
\begin{aligned}
\frac{d\vartheta_4}{d\tau}&=
\frac{\ri}{\pi}\bigg\{\eta-\frac{\pi^2}{12}\,
\big(\vartheta_2^4+\vartheta_3^4
\big)\bigg\}\vartheta_4\,,\\
\frac{d\eta}{d\tau}&=\frac{\ri}{\pi}\bigg\{2\,\eta^2-
\frac{\pi^4}{12^2}
\big(\vartheta_2^8+\vartheta_3^8+\vartheta_4^8 \big) \bigg\}\,.
\end{aligned}
\end{equation}
\end{lemma}
%The system \eqref{var} constitutes a 4-dimensional extension of the
%well-known 3-dimensional Weierstrass--Halphen system on Weierstrass'
%invariants $\g2$, $\g3$, $\eta$ \cite{tannery}:
%\begin{equation}\label{g2g3}
%\frac{d \g2}{d\tau} = \frac{\ri}{\pi}\,
%\Big(8\,\g2\,\eta-12\,\g3\Big)\,,\qquad \frac{d
%\g3}{d\tau}= \frac{\ri}{\pi}\!
%\left(12\,\g3\,\eta-\Mfrac23\,\g2^2\right),\qquad
%\frac{d\eta}{d\tau}=
%\frac{\ri}{\pi}\!\left(2\,\eta^2-\Mfrac16\,\g2\right).
%\end{equation}

In order to obtain the general Abelian differential $R(x,y)\,dx$ on
a PH-curve it is convenient to have a representation for basic
differentials $dy_{\sss\textsc{Pic}}^{\mathstrut}=\dot
y_{\sss\textsc{Pic}}^{\mathstrut}\,d\tau$ and
$dy_{\sss\textsc{Hit}}^{\mathstrut}=\dot
y_{\sss\textsc{Hit}}^{\mathstrut}\,d\tau$ independently of algebraic
form  \eqref{F}.

\begin{theorem}\label{T4}
Let parameters $A,B$ defined by \eqref{nmN} determine the algebraic
solutions \eqref{algPic} and \eqref{param}. Then PH-differentials
$dy=\dot y\, d\tau$, as differentials on corresponding curves
\eqref{F}, are given by the following expressions\/$:$
\begin{alignat}{2}
\dot y&=\pi\ri\vartheta_3^6\,\theta_1^2\,
\frac{y\,(y-1)\,(y-x)}{\vartheta_3^2\,\theta_1^2\,y+
\vartheta_4^2\,\theta_2^2}+
\frac{\pi}{2\ri}\,\vartheta_3^4\,(y-x)^2+\frac{\pi}{2\ri}\,
\vartheta_2^4\,x&&\qquad(\text{Hitchin's curves})\,,\notag\\
\label{diff}\ds\dot y&=\ri \vartheta_2^2
\left(\Dtheta+\pi\ri\Mfrac{\nu}{N}\,\theta_1\right)
\frac{\theta_3\,\theta_4}{\theta_2\,\theta_1^2}^{\strut}\,y+
\pi\ri\vartheta_3^4\,y\,(y-1) &&\qquad(\text{Picard's curves})\,.
\end{alignat}
\end{theorem}

\begin{pf}
This is in fact the Okamoto transformations \eqref{ok} resolved with
respect to $y_x^{}$ which is proportional to $\dot y$. At first
glance Corollary~\ref{C1} contradicts to parametrization
\eqref{param}; whence it follows that function $\Dtheta$ must be a
rational function of $x,y$. With use of theorems of
addition/multiplication for $\zeta$ and the fact that
$\zeta\sim\Dtheta/\theta$ we obtain that
$\Dtheta\big(\frac{\nu}{N}\tau+\frac{\mu}{N}\big|\tau\big)$ is
indeed the rational function of $\theta$'s\footnote{An example:
$\Dtheta\big(\frac14\big)=\frac{\pi}{2}\,(\vartheta_3^2+
\vartheta_4^2) \,\theta_1\big(\frac14\big)$.} over
$\mathbb{C}(\vartheta)$. In other words, Okamoto's transformations
generate through the object $\Dtheta$ the basic Abelian
differentials $dy(\tau)$ independently of $N$. Hitchin's
differentials  $\dot y$ can also be considered as a
$\theta$-function $\tau$-version of the transformations themselves.
When simplifying to form \eqref{diff} we used the equality $\dot
x=-\pi\ri x\,\vartheta_2^4$ being a consequence of Lemma~\ref{L2}.
\hfill $\blacksquare$
\end{pf}

\begin{remark}[\sf Example]\label{R3}
One can consider \eqref{diff} as a Riccati equation $\dot y=
f(\tau)\,y^2+g(\tau)\,y+h(\tau)$ and derive  linear `integrable'
equation of the form $\psi_{\tau\tau}=U(\tau)\,\psi$. By this way we
produce a large (new) family of 2nd order (solvable) \odes\ with
modular coefficients.
\end{remark}

Previous results were obtained with use of some particular
properties of \eqref{p6}, \eqref{param}--\eqref{ok}, and curves;
instead, we can use at once the differential properties of
$\theta$-functions involved in uniformizations \eqref{param}. Below
is a complete description of such properties. It will follow that
differential analysis on PH-curves is just a calculus of the
$\vartheta$- and $\theta$-constants.

\begin{theorem}\label{T5}
Let $A$, $B$ be arbitrary quantities independent of $\tau$. Then the
five theta-constants $\Dtheta,\theta_k(A\tau+B|\tau)$ are
differentially closed over $\mathbb{C}(\vartheta^2,\eta;\pi,A)$\/$:$
\begin{equation}\label{Dtau}
\begin{aligned}
\frac{d\theta_k}{d\tau}
&=\frac{-\ri}{4\pi}\,(\Dtheta+4\,\pi\ri A\,\theta_1)\,
\Dtheta\,\frac{\theta_k}{\theta_1^2}+
\frac{\ri}{2}\,\vartheta_k^2{\cdot}(\Dtheta+
2\,\pi\ri A\,\theta_1)\,
\frac{\theta_\nu\,\theta_\mu}{\theta_1^2}
\\
&\== +\frac{\pi\ri}{4}\,\Big\{
\vartheta_3^2\,\vartheta_4^2{\cdot}\theta_2^2-
\vartheta_k^2\,\vartheta_\mu^2{\cdot}
\theta_\nu^2-
\vartheta_k^2\,\vartheta_\nu^2{\cdot}
\theta_\mu^2\Big\}\,\frac{\theta_k}{\theta_1^2}+
\frac{\ri}{\pi}\bigg\{
\eta+\frac{\pi^2}{12} \big(\vartheta_3^4+
\vartheta_4^4\big)\bigg\}{\cdot} \theta_k\,,
\\
\frac{d\Dtheta}{d\tau}&=\frac{\ri}{4\pi}\,
(3\,\Dtheta+4\,\pi\ri A\,\theta_1) \bigg\{
\pi^2\,\vartheta_3^2\,\vartheta_4^2{\cdot}
\frac{\theta_2^2}{\theta_1^2}+
4\,\eta+\frac{\pi^2}{3}\,\big(\vartheta_3^4+\vartheta_4^4\big)\bigg\}
\\
&\== -\frac{\ri}{4\,\pi}\,(\Dtheta+4\,\pi\ri A\,\theta_1)\,
\frac{(\Dtheta)^2}{\theta_1^2} -\frac{\ri}{2}\,\pi^2
\vartheta_2^2\,\vartheta_3^2\,\vartheta_4^2{\cdot}
\frac{\theta_2\,\theta_3\,\theta_4}{\theta_1^2}\,,
\end{aligned}
\end{equation}
where $k=1,2,3,4$ and symbols $\nu$, $\mu$ stand for
\begin{equation}\label{munu}
\nu=\frac{8\,k-28}{3\,k-10}\,,\qquad\mu=\frac{10\,k-28}{3\,k-8}\,.
\end{equation}
\end{theorem}

\begin{pf}
Denote by prime $'$ the derivative with respect to $z$. Then the
logarithmic derivatives ${\ln'}\theta_k(z|\tau)$ appear in quadratic
identities of the type
$\Dtheta[\alpha]\,\theta_\beta-\Dtheta[\beta]\,\theta_\alpha\sim
\theta\,\theta$ which are sometimes present in the old literature
\cite{baruch,tannery,weber}. In order to write them in general form
we use notation \eqref{munu} which tells us that indices
$(k,\nu,\mu)$ run over cyclic permutations of the series $(2,3,4)$
when $k=2,3,4$. The $\theta,\theta'$-identities above then can be
rewritten as follows
$$
\frac{\Dtheta[k]}{\theta_k}-\frac{\Dtheta}{\theta_1}=
-\pi\vartheta_k^2\cdot\frac{\theta_\nu\,\theta_\mu}
{\theta_k\,\theta_1}\,.
$$
This gives derivatives $\Dtheta[2,3,4]$ in terms of $\theta$,
$\Dtheta$ and the formula works also well under $k=1$  since
$\vartheta_k=\vartheta_1\equiv0$ and the terms $\sim\vartheta_k$
drop out. Converting the obvious elliptic identity
$(\sigma')'=(\sigma')^2/\sigma-\sigma\,\wp$ into $\theta$-functions,
we obtain the expression for $(\Dtheta)'$ through $\theta$ and
$\Dtheta$ itself:
$$
\frac{\partial\Dtheta}{\partial z}=
\frac{(\Dtheta)^2}{\theta_1}-\pi^2\vartheta_3^2\,\vartheta_4^2
\cdot \frac{\theta_2^2}{\theta_1}^{\strut}-
4\,\bigg\{\eta+\frac{\pi^2}{12}\big(\vartheta_3^4+\vartheta_4^4\big)
\bigg\} \cdot\theta_1\,.
$$
Invoking the heat equations
$4\pi\ri\theta_\tau=\theta_{\mathit{zz}}$,
$4\pi\ri\Dtheta[\tau]=\Dtheta[\mathit{zz}]$ and making the change
$z\mapsto A\,\tau+B$, we arrive, upon simplification, at equations
\eqref{Dtau}.  Lemma~\ref{L2} provides their differential closure
and serves also the case $A\tau+B=0$ as a limiting case of
Eqs.~\eqref{Dtau}.\hfill $\blacksquare$
\end{pf}

The functions $\theta_k(A\,\tau+B|\tau)$ are the continual
generalizations of $\theta$'s with discrete characteristics and this
theorem is  a particular case of more general differential
properties of $\theta$-functions which are exhaustively described in
\cite{br4}. An addition of fourth equation in \eqref{var} is an
important point. In a particular case, when $A$ and $B$ are chosen
to be \eqref{nmN},  we could formally get by without second equation
in \eqref{Dtau} but form of the first one would then depend on $N$.
In other words, the Okamoto transformations \eqref{ok}, Abelian
differentials \eqref{diff}, curves themselves \eqref{F}, and
Fuchsian equations \eqref{Qxy}, being rewritten in the language of
$\theta$-functions, constitute numerous and rather sophisticated
$\eta,\vartheta,\theta$-identities and their differential
consequences.

\begin{remark}\label{R4}
It is a good exercise  to check the original solution of Hitchin
\eqref{morazm} with direct use of Theorem~\ref{T5} and
Lemma~\ref{L2}.
\end{remark}

\section{Group transformations\label{groups}}

\subsection{Basic $\theta$-transformations\label{groupstheta}}

In order to obtain the group properties of the appeared automorphic
functions we need transformation properties of their ingredients,
\ie, $\theta,\theta'$-constants. The transformation for function
$\theta_1$ is known \cite{tannery,weber} and usually written in the
form
\begin{equation}\label{ab11}
\theta_1\Big(
\Mfrac{z}{c\,\tau+d}\Big|\Mfrac{a\,\tau+b}{c\,\tau+d}\Big)=
\Aleph\cdot\!\sqrt{c\,\tau+d\:}\, \,\re_{\strut}^{\frac{\pi\ri
c\, z^2}{c\,\tau+d}}\,\theta_1(z|\tau)\,,
\end{equation}
where $\Aleph$ denotes some eighth root of unity and
$\big(\begin{smallmatrix}a&b\\c&d
\end{smallmatrix}\big)\in
\bo\Gamma(1)$. However we shall require more detailed information.
Let $[p]$ stand for integer part of the number $p$ and suppose that
$c$ is normalized to be positive: $c>0$.

\begin{theorem}[\cite{br4}]\label{T6} The
$\bo\Gamma(1)$-transformation law of the general $\theta$-function.
Let $\theta\AB{\alpha}{\beta}$ be the theta-series with integer
characteristics \eqref{theta} and $n\in\mathbb{Z}$. Then
\begin{align}\notag
\theta\AB{\alpha\sm1}{\beta} (z|\tau+n)&=
\ri_{\mathstrut}^{\frac n2\,(1\sm\alpha^2)}
\cdot \theta\AB{\alpha-1}{\beta+n\alpha} (z|\tau)\,,\\
\label{ab}\ds
\theta\Larger{\AB{\tilde\alpha\sm1}{\tilde\beta\sm1}}
\Big(\Mfrac{z}{c\,\tau+d}\Big|
\Mfrac{a\,\tau+b}{c\,\tau+d}\Big)&=
\bo{\mathfrak{E}}_{\alpha\beta}\:\Aleph\cdot\!
\sqrt{c\,\tau+d\,}
\re^{\frac{\pi\ri cz^2}{c\,\tau+d}}_{\mathstrut}\,
\theta\AB{\alpha\sm1}{\beta\sm1}(z|\tau)\,,
\end{align}
where multipliers $\Aleph$ and $\bo{\mathfrak{E}}_{\alpha\beta}$
read
\begin{align}
\bo{\mathfrak{E}}_{\alpha\beta}&=
\ri^{\frac12\mbig[0]\{
2\,\alpha\,(\beta\,b\,c-d+1)-\beta\,c\,(\beta\,a-2)-
\alpha^2\,d\,b\mbig[0]\}}_{\mathstrut}\,,\notag\\
\label{aleph}
\Aleph&= \exp3\pi\ri  \mbig[11]\{\frac{a-d}{12\,c}
-\frac{d}{6}(2c-3)+\frac{c-1}{4}\mathfrak{sign}
(d)-\frac14+\frac1c\cdot
\SUM{k}{\raisebox{-0.07em}{\smaller\textup\textbar}\!\!
{[c\mspace{-0.9mu}/\mspace{-2mu}d]\mspace{-2.5mu}
\raisebox{-0.07em}{\smaller\textup\textbar}\mspace{-1mu}
+\mspace{-1mu}1}}{{\scriptstyle c-1}}
\Big[\Mfrac{d}{c}k\Big]\,k \mbig[11]\}\,,
\end{align}
and characteristics $(\alpha,\beta)$, $(\tilde\alpha,\tilde\beta)$
are related through the linear transformations
\begin{equation}\label{abAB} \Bigg\{\;
\begin{alignedat}{6}
\widetilde\alpha&= &&d\,&\alpha &{}-{}c\,&\beta\,,\quad\qquad
\alpha&{}=a\,&\widetilde\alpha&{}+c\,&\widetilde\beta\\
\widetilde\beta &=-&&b\,&\alpha &{}+{}a\,&\beta\,,\quad\qquad
\beta&{}=b\,&\widetilde\alpha&{}+d\,&\widetilde\beta
\end{alignedat}\;
\Bigg\}.
\end{equation}
\end{theorem}

The main point here is an explicit and self-contained form of
$\theta$-transformation. The multiplier ${\Aleph}$ is a common
quantity to all the $\theta$'s since it does not depend on
characteristics $(\alpha,\beta)$, whereas $\bo{\mathfrak{E}}$ does.
Classical uniformizing functions are determined by classical
$\vartheta$-constants and the object
$\bo{\mathfrak{E}}_{\alpha\beta}(a,b,c,d)$. In turn, the well-known
property of $\theta$-characteristics
\begin{equation}\label{a2m}
\theta\AB{\alpha+2m}{\beta+2n}(z|\tau)=
(-1)^{n\alpha}\,\theta\AB{\alpha}{\beta}(z|\tau)
\end{equation}
define congruence properties of monodromies. Amongst the curves
which have appeared, the three such functions encounter, each is
universal uniformizing one:\footnote{In particular, third of these
functions uniformizes the simplest spectral curve $w^3=v^5-\,v$ of a
tetrahedrally symmetric non-singular monopole of charge 3
\cite[Theorem~1]{hitchin2}. See \cite[p.~258]{br} for
parametrization and \cite{braden} for special discussion of this
monopole.}
\begin{equation}\label{xuv}
x=\frac{\vartheta_4^4(\tau)}{\vartheta_3^4(\tau)}\,,\qquad
u=\frac{\vartheta_4^2(\tau)}{\vartheta_3^2(\tau)}\,,\qquad
v=\frac{\vartheta_4(\tau)}{\vartheta_3(\tau)}\,.
\end{equation}
We have, however, general $\theta',\theta(A\,\tau+B|\tau)$-constants
under restriction \eqref{nmN} and therefore we should obtain
automorphy factors for this kind $\theta,\theta'$-constants. If one
considers Hitchin's algebraic solutions we should use the
transformation law for $\Dtheta$ which is derived from \eqref{ab11}:
$$
\Dtheta\!\Big(
\Mfrac{z}{c\,\tau+d}\Big|\Mfrac{a\,\tau+b}{c\,\tau+d}\Big)=
\Aleph\cdot\sqrt{c\,\tau+d\:}\,\re_{\strut}^{\frac{\pi\ri
c z^2}{c\,\tau+d}}\,\mbig\{(c\,\tau+d)\,\Dtheta(z|\tau)+
2\,\pi\ri c\,\,z\,\theta_1(z|\tau) \mbig\}\,.
$$
From the transformations above  it follows immediately that the set
of new $\vartheta,\theta$-constants \eqref{nmN} is closed with
respect to $\bo\Gamma(1)$ under fix $N$:
\begin{equation}\label{Gamma1tau}
\begin{aligned}
\theta_{\alpha\beta}\Big(\Mfrac{\nu}{N}\tau+\Mfrac{\mu}{N}\Big|\,\tau
\Big)\! \quad\overset{\widehat{\bo\Gamma(1)}(\tau)}
{{\scalebox{3.2}[1]{$\mapsto$}}} \quad\!&
\theta_{\alpha\beta}\Big(\Mfrac{\nu}{N}\,\Mfrac{a\,\tau+b}{c\,\tau+d}
+\Mfrac{\mu}{N}\Big| \Mfrac{a\,\tau+b}{c\,\tau+d}
\Big)\\
\cong\:& \theta_{\alpha'\beta'}\Big( \Mfrac{\nu
a+\mu c}{N}\,\tau+\Mfrac{\nu b+\mu d}{N}\Big|\,\tau \Big)=
\theta_{\alpha'\beta'}
\Big( \Mfrac{\nu'}{N}\,\tau+\Mfrac{\mu'}{N}\Big|\,\tau \Big)\,,
\end{aligned}
\end{equation}
so we can build transformations for $\theta$-quotients.

In order to obtain  the congruence for automorphism group of a
function we should require that $\vartheta,\theta$-ratio defining
the function transforms into itself; to do this, formulae
\eqref{ab}, \eqref{abAB}, \eqref{a2m}, and multiplier
$\bo{\mathfrak{E}}$ provide all the required information.
Furthermore, all the automorphic functions arisen from PH-curves are
algebraically related to the function $x(\tau)$. Hence their
automorphisms  are commensurable with $\bo\Gamma(2)$ in
$\bo\Gamma(1)$. This simplifies the analysis because transformation
\eqref{ab}, when
\begin{equation}\label{abcd2}
\Big(\,\begin{matrix}a&b\\c&d
\end{matrix}\,\Big)=
\mbig[7](\,\begin{matrix}2n+1&2m\\2p&2q+1
\end{matrix}\,\mbig[7])\in\bo\Gamma(2)\,,\qquad
n,m,p,q\in\mathbb{Z}
\end{equation}
splits into the separate formulae for all the
$\theta_k$:
\begin{alignat*}{2}
\theta_1\!\Big(\Mfrac{z}{c\,\tau+d}\Big|
\Mfrac{a\,\tau+b}{c\,\tau+d}\Big)&{}={}&\Aleph \cdot\!
\sqrt{c\,\tau+d\,}\,\, \re^{\!\frac{\pi\ri c\,z^2}{c\,\tau+d}}\,
\theta_1(z|\tau)\,,\\
\theta_2\!\Big(\Mfrac{z}{c\,\tau+d}\Big|
\Mfrac{a\,\tau+b}{c\,\tau+d}\Big)&{}={}&\ri^
{2q(p-1)+p}\,\Aleph\cdot\! \sqrt{c\,\tau+d\,}\,\,
\re^{\!\frac{\pi\ri c\,z^2}{c\,\tau+d}}\, \theta_2(z|\tau)\,,
\\
\theta_3\!\Big(\Mfrac{z}{c\,\tau+d}\Big|
\Mfrac{a\,\tau+b}{c\,\tau+d}\Big)&{}={}& \ri^{2q(p+1)-m(2n+1)+p}\,
\Aleph\cdot\! \sqrt{c\,\tau+d\,}\,\, \re^{\!\frac{\pi\ri
c\,z^2}{c\,\tau+d}}\, \theta_3(z|\tau)\,,
\\
\theta_4\!\Big(\Mfrac{z}{c\,\tau+d}\Big|
\Mfrac{a\,\tau+b}{c\,\tau+d}\Big)&{}={}&
\ri^{2n(m-1)-m}\,\Aleph\cdot\! \sqrt{c\,\tau+d\,}\,\,
\re^{\!\frac{\pi\ri c\,z^2}{c\,\tau+d}}\, \theta_4(z|\tau)\,,
\end{alignat*}
where $\Aleph$ is recomputed according to \eqref{aleph} and
\eqref{abcd2}. These and many other details and useful properties of
$\theta$-functions and $\vartheta$-constants can be found in
\cite{br4}.

\subsection{Congruences for Picard's groups\label{groupsPicard}}

The formulae above, supplemented with the rule \eqref{a2m}, allow us
to derive the congruence properties of uniformizing functions under
question. We shall restrict our consideration to Picard's curves.

\begin{theorem}\label{T7}
The automorphism group $\G_y$ of Picard's uniformizing function
\begin{equation}\label{Pic}
y(\tau)=-\frac{\vartheta_4^2(\tau)}{\vartheta_3^2(\tau)}\,
\frac{\theta_2^2\big(\frac{\nu\,\tau+\mu}{2\,N}\big|\tau \big)}
{\theta_1^2\big(\frac{\nu\,\tau+\mu}{2\,N}\big|\tau \big)}
\end{equation}
coincides with the group $\G^{\s{(N)}}_{\sss\textsc{Pic}}$
uniformizing the Picard curve and is a free congruence subgroup of
$\bo\Gamma(2)$ with one of the following representations\/$:$
\begin{equation}\label{nmpq}
\G^{\s{(\bo{N})}}_{\sss\textsc{Pic}}\colon\qquad\G=
\mbig[7](
\begin{matrix}2\bo{N}n+1&2\,m\\2\bo{N}p&2\,
\bo{N}q+1
\end{matrix} \mbig[7])\,,\qquad
\G^{\s{\bo\top}}=\mbig[7](
\begin{matrix}2\bo{N}n+1&2\bo{N}p\\2\,m&2\,
\bo{N}q+1
\end{matrix} \mbig[7])\,,
\end{equation}
\begin{equation}\label{cong}
2\,m\,p-2\,\bo{N}q\,n=q+n\,,\qquad n,m,p,q\in\mathbb{Z}\,,
\end{equation}
where $(\nu,\mu,N)\ne (0,\pm 1,\pm2)$ and $\bo{N}\leqslant N$ is
computed through $(\nu,\mu,N)$. The group
$\G^{\s{(N)}}_{\sss\textsc{Pic}}$ has the same topological genus as
genus of Picard's curve $F_{\!\s{(\bo N)}}(x,y)=0$ under
$(\nu,\mu)=(0,1)$. Rank and generators of the group
$\G^{\s{(N)}}_{\sss\textsc{Pic}}$ are computable.
\end{theorem}

\begin{pf}
The case $N=1$ is trivial and we assume that $N>1$. Since
$\vartheta$- and $\theta$-constants transform separately into
themselves, we obtain, based on $\bo\Gamma(2)$-split
$\theta$-transformations above and invariancy condition for the
theta-ratios in \eqref{Pic}, that group $\G_y$ is not only
commensurable with $\bo\Gamma(2)$ but is its subgroup:
$\G_y\in\bo\Gamma(2)$. We know also that $\bo\Gamma(2)=\G_x$ and
$(x,y)$ are generators of the function field on the curve
$F_{\!{\s{(N)}}}(x,y)=0$; hence it follows that the image of this
curve in $\Hp$ coincides with geometrical polygon for the group
$\G^{\s{(N)}}_{{\sss\textsc{Pic}}}$. Therefore
$\G^{\s{(N)}}_{\sss\textsc{Pic}}= \G_y\cap\G_x=\G_y$.

There is an exceptional case of $z$ under which
$\theta_2(z|\tau)=\pm\theta_1(z|\tau)$. The known property
$\theta_1\big(z+\frac12\big|\tau\big)=\theta_2(z|\tau)$ implies that
this is the case
$\theta_2\big({\pm}\frac14\big)=\pm\theta_1\big({\pm}\frac14\big)$
and therefore $(\nu,\mu,N)= (0,\pm 1,\pm2)$ (see {\sf
Example~\ref{E2}} further below). There are two sets of parameters
$(\nu,\mu,N)$ in \eqref{Pic}:
\begin{alignat}{2}\label{case1}
&\theta\Big(\Mfrac{\nu\,\tau+\mu}{\bo{N}}
\Big)\,,\qquad&&\bo{N}=2,3,4,\...
\\\intertext{and}\label{case2}
&\theta\Big(2^k\,\Mfrac{\nu\,\tau+\mu}{\bo{N}}
\Big)\,,\qquad&&\bo{N}=3,5,7,\...\,,\quad k=1,2,3,\...\,,
\end{alignat}
where $\nu$ and $\mu$ are not even simultaneously. Both of the cases
\eqref{case1}, \eqref{case2} come from  Weierstrass'
$\wp\mbig(\frac{\nu\,\tau+\mu}{N}\big|\tau \mbig)$ or
$\wp\mbig(2^k\,\frac{\nu\,\tau+\mu}{N}\big|\tau \mbig)$ with the same
meaning of $\nu,\mu$. Let $R_N$ denote a rational function of $\wp$
determining the multiplication theorem for $\wp(Nz)$. Therefore
\begin{equation}\label{r1}
R_N\Big( \wp\Big(\Mfrac{\nu\,\tau+\mu}{N}
\mbig|\tau\Big)\Big)=\wp(\nu\,\tau+\mu|\tau)=\cdots
\end{equation}
and this expression may have only one of the three values
\begin{equation}\label{r2}
\cdots=\mbig\{\wp(1|\tau),\;\wp(\tau|\tau),\;
\wp(\tau+1|\tau)\mbig\}\,,
\end{equation}
so that all the values $\nu,\mu$ are equivalent to the three cases
$\nu\,\tau+\mu=\{1,\tau,\tau+1\}$:
$$
\wp(1|\tau)=\frac{\pi^2}{12}\,\vartheta_3^4(\tau)\,(x+1)\,,
\qquad
\wp(\tau|\tau)=\frac{\pi^2}{12}\,\vartheta_3^4(\tau)\,(x-2)\,,
\qquad
\wp(\tau+1|\tau)=\frac{\pi^2}{12}\,\vartheta_3^4(\tau)\,(1-2\,x)\,.
$$
(Hence, in order to get all PH-curves one needs to use only
multiplication theorems\footnote{No need to involve addition
formulae as mentioned in \cite[Lemma~3]{mazz}; multiplication
formulae have nice and effective recurrences \cite[{\bf
III:}~p.~105]{tannery} and are much more effective in use than the
addition ones.} for $\wp$). The two latter cases produce
automorphisms $\G_y$ which are conjugated to the case $\wp(1|\tau)$
since $\wp(1|\tau)\mapsto\wp(\tau|\tau)$ under
$\tau\mapsto-\frac1\tau$ and $\wp(1|\tau)\mapsto\wp(\tau+1|\tau)$
under $\tau\mapsto\frac{\tau}{\tau+1}$. Furthermore, in this case
$\wp\mbig(2^k\,\frac{\nu\,\tau+\mu}{\bo{N}}\big|\tau
\mbig)=R_{2^k}\!\mbig[3](\wp\big(\frac{\nu\,\tau+\mu}{\bo{N}}
\big|\tau \big)\mbig[3])$ and  we get a rational function of the
algebraic one corresponding to the case
$\wp\big(\frac{\nu\,\tau+\mu}{\bo{N}}\big|\tau \big)$ with $\nu,\mu$
equal to 0 or 1. Therefore genera of groups
$\G^{\s{(N)}}_{\sss\textsc{Pic}}$ and $\G^{\s{(\bo
N)}}_{\sss\textsc{Pic}}$ coincide. Again, all the groups will be
conjugated to the group with $(\nu,\mu)=(0,1)$ and topological genus
of  polygon for $\G^{\s{(N)}}_{\sss\textsc{Pic}}$ coincides with the
genus of the curve \mbox{$F_{\!\s{(\bo N)}}(x,y)=0$}.

Transformations \eqref{Gamma1tau} imply the following congruence
conditions on entries of \eqref{abcd2}:
$$
2^k\,(\nu\,n+\mu\,p)=\bo{N}P\,,\qquad
2^k\,(\nu\,m+\mu\,q)=\bo{N}Q\,,\qquad P,Q\in\mathbb{Z}\,.
$$
By previous argument, we may consider only the case
$\nu\,\tau+\mu=1$, where $\bo{N}$ is a free  parameter being an odd
number or $\bo{N}\in\mathbb{Z}$ under $k=0$. Therefore $P,Q\sim 2^k$
and we have $p,q\sim \bo{N}$.  The unimodular condition $\det\G=1$
yields $2\bo{N}(m\,p-n\,q)=\bo{N}q+n$ and therefore $n\equiv 0\!\mod
\bo{N}$. Replacing $n$ with $\bo{N}n$, we arrive at \eqref{cong} and
the first set $\G$ of matrices in \eqref{nmpq}. Above-mentioned
transformations of $\wp$'s have been formed by the two ones:
$$
\tau\mapsto\tau+1\,,\qquad\tau\mapsto-\frac1\tau\,.
$$
Conjugating $\G$ by the first of these transformations, we get $\G$
itself and the second one produces$(\G^{\s{\bo\top}})^{\sm1}
=\G^{\s{\bo\top}}$; it is the second matrix in \eqref{nmpq}.

There is a general algorithm concerning  subgroups of free groups
and it is known as the famous Reidemeister--Schreier rewriting
process \cite{solitar,lyndon}. As applied to our automorphisms, the
algorithm is somewhat too general since it deals with abstract group
presentations meanwhile we are concerned  with the matrix
monodromies. Since $\G^{\s{(N)}}_{\sss\textsc{Pic}}$ is a subgroup
of the free group $\bo\Gamma(2)$, the monodromy $\G_y$ (in order to
be free) has to have a parabolic element $T_0$ and this generator,
being a generator of the \emph{global} monodromy $\G_y$, must
correspond to punctures at $\Hp\!\big/\G_y$ and $\Hp\!\big/\G_x$.
Obviously, $T_0$ is a power of some parabolic element from
$\bo\Gamma(2)$. We thus get a polygon for $\G_y$ as a set of copies
of $x$-quadrangles. Generators of the global monodromy $\G_y$ are
easily obtainable by geometric analysis of this `big' polygon
supplemented with solution of congruence \eqref{cong}.\hfill
$\blacksquare$
\end{pf}

The following table illustrates genera and degrees of the covers
$x\mapsto y$ for several PH-curves; we took the maximal $y$-degree
connected component of solutions determined by formulae
\eqref{r1}--\eqref{r2}:
\begin{equation}\label{genera}
\begin{array}{r|c|c|c|c|c|c|c|c|c|c|c|c|c}
N&2&3&4&5&6&7&8&9&10&11&12&13&\cdots\\\hline
\text{\small genus}&0\vbox to2.5ex{}&0&0&1&1&4&5&7&9&16&
\textit{\ttt{13}}&25&\cdots\\\hline
\text{\small number of $y$-sheets}
&2\vbox to2.5ex{}&4&8&12&16&24&32&36&48&60&64&84&\cdots
\end{array}
\end{equation}
For prime $N$'s these genera fit in the Barth--Michel formula
$g=\frac14(N-3)^2$ \cite{barth} and Hitchin did show that the
formula works for algebraic solutions defined by Poncelet's polygons
\cite{hitchin3}. The general formula for \eqref{genera} is unknown;
the table has some `sporadic' entries like 13 (italicized  in
\eqref{genera}).

It remains to consider the exceptional case in Theorem~\ref{T7}. We
shall do this more fully since this case demonstrates the way of
getting formulae.

\begin{example}\label{E2}
Let us consider function $u(\tau)$ in \eqref{xuv}. Imposing the
invariancy condition for this ratio, we obtain equations on entries
of the transformation:
\begin{equation}\label{p02}
\ri^{2p}=1\quad \hence\quad p=\{0,2\}\,,\qquad a\,d-b\,c=1\quad
\hence\quad 2\,(m\,p-q\,n)=q+n\,.
\end{equation}
Integral solutions to these equations completely determine the group
$\G_u=\text{Aut}\,u(\tau)$, that is the monodromy group of   Heun's
equation for Legendre's modulus $u=k(\tau)$:
\begin{equation}\label{lemn}
\Psi''=-\frac14\,\frac{(u^2+1)^2}{u^2(u^2-1)^2} \,\Psi\,.
\end{equation}
It is obvious that we have to have three generators for this  group
and they are obtainable  from \eqref{p02}. The group $\G_u$ is a
freely generated one  with index 2 in $\bo\Gamma(2)$ and we need
only find any three solutions of \eqref{p02}  being integers $n,m,q$
nearest to zero.  The pairing of neighboring quadrangles for
$\bo\Gamma(2)$ implies that we should choose only parabolic
representatives. One easily computes:
\begin{gather*}
\begin{aligned}
\G_u &=\mbig\langle T_{(1)},\,T_{(0)},\,T_{(\infty)}\mbig\rangle=
\mbig\langle U,\,S^2,\, (SU^{\sm1})^2\mbig\rangle\\
&= \bigg\langle
\mbig[7](\,\begin{matrix}1&\!\!2\\0&\!\!1\end{matrix}
\,\mbig[7]),\quad
\mbig[7](\,\begin{matrix}1&\!\!0\\4&\!\!1\end{matrix}
\,\mbig[7]),\quad
\mbig[7](\,\begin{matrix}3&\!\!-4\\4&\!\!-5\end{matrix}
\,\mbig[7])\bigg\rangle\,,
\end{aligned}\\
u(\ri \infty)=1\,,\quad u(0)=0\,,\quad u(1)=\infty\,,\quad
u\mbig[5](\mfrac12\mbig[5])=-1\,,
\end{gather*}
where we have represented $T$'s in form of  decompositions by the
standard $\bo\Gamma(2)$-generators:
$U=\big(\begin{smallmatrix}1&2\\0&1
\end{smallmatrix}\big)$ and $S=\big(\begin{smallmatrix}1&0\\2&1
\end{smallmatrix}\big)$. This matrix
representation for monodromy of \eqref{lemn} has been attached to a
certain basis of solutions $(\Psi_1, \Psi_2)$. Such a basis is
calculated, as usual in uniformization,  through the uniformizing
Hauptmodul \cite{ford}: $\Psi\sim\sqrt{\dot x}$. We then get
\begin{equation}\label{uK}
\sqrt{u_\tau}=\sqrt{\frac{\pi}{2\ri}\,
\frac{\vartheta_4^2}{\vartheta_3^2}\,\vartheta_2^4}\,\sim\Psi_1=
\sqrt{u\,(u^2-1)}\,\ellK'(u),\qquad \Psi_2=\ri
\sqrt{u\,(u^2-1)}\,\ellK(u)\,.
\end{equation}
Matrix monodromy of this function pair \eqref{ma} is exactly group
$\G_u$ above.
\end{example}

\section{Towers of uniformizable curves. Examples\label{towers}}

We now pass to consequences. Having two (not necessarily PH's)
curves in hand, $F_1(x,y)=0$ and $F_2(x,\tilde y)=0$, we may
formally eliminate variable $x$ and get one more curve $F(y,\tilde
y)=0$. In particular, \emph{any algebraic PH-solution is
algebraically related to any other PH-algebraic solution}. Passing
now from algebraic curves  to their uniformizations $y(\tau)$,
$\tilde y(\tau)$, we see that such an elimination is compatible with
parametrizations if and only if the one common function $x(\tau)$
appears in parametrizations of both of the curves, \ie, $x(\tau)$ is
the universal uniformizing function for branch $x$-points of $F_1$
and $F_2$. This seemingly trivial procedure, once transformed into
the `$\tau$-representation', leads to rather nontrivial
consequences. Namely, new and completely uniformizable curves of
higher genera. Indeed, let $x=x(\tau)$ be the universal
uniformization for curves $F_1$ and $F_2$ and functions $y(\tau)$
and $\tilde y(\tau)$ are known. Then we obtain not only
parametrization of the curve $F(y,\tilde y)=0$, but Fuchsian
equations, groups, and differentials as well. The general
explanation of this fact is that automorphism group $\G_z$ of some
modular function $z(\tau)$ has a lot of subgroups of higher genera
even though modular equation defining the function $z(\tau)$ has a
zero genus. In this respect the equation itself is just a particular
case of this infinite family. We thus construct curves of nontrivial
genera without seeking for complicated subgroups of $\bo\Gamma(1)$
having nontrivial genera. For example, Klein's curves
$F\big(J(\tau),J(N\tau)\big)=0$ and their groups $\Gamma_0(N)$ have
zero genera for $N$ up to 10, whereas $\Gamma_0(25)$ has again the
genus zero; it is known that size of these equations comes into
extremely rapid growth. See \cite{maier2} for details and explicit
formulae. As numerous examples show, the curves obtained from
PH-series are simpler than the majority of classical modular
equations. In addition to all this, Theorem~\ref{T5} provides a
complete differential apparatus for these curves.
\begin{example}\label{E3}
The relation between $u(x)=\sqrt{x}$ and Picard's curve $y(x)$ of
level $N=5$ is a curve of genus $g=5$:
\begin{equation}\label{g5}
\begin{aligned}
&(16\,y^2-20\,y+5)\, u^{12}-2\, y\,(40\,y^2-47 \,y+10)\, u^{10}\\
&\quad{}+y^2\,(64\,y^5-240\,y^4+360\,y^3 -105\,y^2-80\,y+16)\, u^8\\
&\quad{}-20\, y^5\,(8\,y^3- 28\,y^2+39\,y-18)\, u^6+
5\,y^6\,(28\,y^3-89 \,y^2+112\,y-48)\,u^4\\
&\quad{}-2\,y^7\,(5\,y-4)\,(5\,y^2-10\,y+8)\,
u^2+y^{12}=0\,.
\end{aligned}
\end{equation}
Uniformizing functions $u(\tau)$ and $y(\tau)$ for this curve read as
follows
$$
u=\frac{\vartheta_4^2(\tau)}{\vartheta_3^2(\tau)}\,, \qquad
y=-\frac{\vartheta_4^2(\tau)}{\vartheta_3^2(\tau)}\,
\frac{\theta_2^2\big(\frac{1}{10}
\big|\tau\big)}{\theta_1^2\big(\frac{1}{10}\big|\tau \big)}
$$
and have monodromies of genus zero and unity respectively (see table
\eqref{genera}). The curve \eqref{g5} is not hyperelliptic but can
be realized as a cover of a torus. Indeed, equation \eqref{g5} has
the obvious sheet interchange symmetry $u\mapsto -u$ and therefore
we may consider \eqref{g5} as a cover of the plane $u^2$, which is
the $x$-plane again. A simple calculation shows that this cover is
an arithmetic torus equivalent to the form $w^2=z^3-12\,z-11$ having
Klein's $J$-invariant equal to $\frac{256}{135}$. Incidentally it
should be remarked that the above-written PH-parametrization of the
torus is highly elementary, whereas its standard $\wp$-Weierstrass'
one is too tremendous to display here.
\end{example}

\subsection{Non-3-branch covers\label{Hyper}}

As well as being a very wide class of completely describable curves,
PH-curves provide a large number of curves of not so special form as
themselves. Curves of this family are also of interest because
these, contrary to Painlev\'e curves, have greater than three branch
points and, therefore, do not belong to the well-known class of the
Bely\u\i\ curves. Formula \eqref{g5} is a good example but even
simpler PH-curves yield a rich theory with nice consequences.

Let us consider two simplest Picard's cases corresponding to
$N=\{2,3\}$ (by virtue of Corollary~\ref{C1} we could equally well
take Hitchin's formulae). They produce two solutions $u(x)$ and
$y(x)$ defined by the curves (Hitchin--Dubrovin)
\begin{equation}\label{uy}
u^2=x\,,\qquad y^4-6\,x\,y^2+4\,x\,(x+1)\,y-3\,x^2=0\,;
\end{equation}
second of these solutions corresponds to $A\,\tau+B=\frac16$.
Solutions \eqref{uy} lie on the elliptic curve
\begin{equation}\label{tor-yu}
y^4-u^2\,(6\,y-4)\,y+u^4\,(4\,y-3)=0
\end{equation}
although both $u(\tau)$ and $y(\tau)$ have zero genus monodromies
(see table \eqref{genera} again). In order to show how general
Fuchsian equations \eqref{Qxy} may look we exemplify an equation for
the function $y(\tau)$. From  Theorem~\ref{T3} we derive
$$
\![5][y,\,\tau]=-\frac18\,{\frac{8\,x\,(y-2)(y^2-9\,y+9)\,y
+16\,y^6+27\,y^5+95\,y^4-415\,y^3+465\,y^2-288\,y+
108}{y^2\,(4\,y-3)(y+3)^2\,(y-1)^3}}\,;
$$
this equation can serve as a nontrivial example of solvable Fuchsian
equation \eqref{fuchs} with algebraic coefficients. Turning back to
the torus \eqref{tor-yu}, we can transform it into any of standard
forms. Such a kind manipulations $(y,u)\rightleftarrows (z,w)$ have
long been algorithmized and we adopt Riemann's form
\begin{equation}\label{TOR}
w^2=z\,(z-1)(z+3)\,.
\end{equation}
It is obtainable from \eqref{tor-yu} by the following birational
isomorphism
\begin{equation}\label{tor}
\begin{aligned}
y&=-\frac{1}{4}\,\frac{w^2}{z}\,,&\qquad\quad
z&=2-3\,\frac{u^2}{y}-\frac{u^2-1}{y-1}\,,\\
u&=\frac14\,\big(1-z^{\sm1}\big)\,w\,,&
w&=\frac{(5\,y-3)\,u^2-y^3-y^2}{u\,(u^2-1)}\,.
\end{aligned}
\end{equation}

That we have a torus in a canonical form does not mean that we shall
arrive at a Fuchsian equation with singularities located precisely
at four branch-places of canonical structures for elliptic tori.
However the simplicity of \eqref{TOR} tells us that function $z$ is
certain to satisfy a simple Fuchsian equation. This is so indeed and
we obtain (Lemma~\ref{L1}) that expression
\begin{equation}\label{z}
z(\tau)=2+3\,\frac{\vartheta_4^2}{\vartheta_3^2}\cdot
\frac{\theta_1^2}{\theta_2^2}-
\frac{\vartheta_2^2}{\vartheta_3^2}\cdot
\frac{\theta_1^2}{\theta_4^2} \,,
\end{equation}
where $\theta\DEF\theta\big(\frac16\big|\tau \big)$, solves an
equation which turns out to be very elegant:
\begin{equation}\label{Qz}
[z,\,\tau]= -\frac12\,\frac {(z^2+3)^4}{(z^5-10\,z^3+9\,z)^2}
\qquad\FED\bcQ(z)\,.
\end{equation}
It is a good exercise to check this equation employing
Theorem~\ref{T5} and Lemma~\ref{L2}; not by Lemma~\ref{L1}.
Eq.~\eqref{z} has regular singular points $E_k=\{0,\pm1, \pm3\}$ and
$E_6=\infty$ so that $z^5-10\,z^3+9\,z$ has $z=E_k$ for roots. Since
$$
\bcQ(z)=-\frac12\,\sum_{k=1}^5\frac{1}{(z-E_k)^2}+
\frac{2\,z^2-6}{(z-1)(z+1)(z-3)(z+3)}\,,
$$
all the six points $\{E_k\}$ correspond to  punctures. Correlating
this equation with torus \eqref{TOR}, we could treat the set
$\{E_k\}$ as the fact that torus \eqref{TOR} has `superfluous'
punctures at points $z=\{{-}1,+3\}$ but more correct explanation, as
a continuation and illustration to Remark~\ref{R1}, reads as
follows.

\begin{remark}\label{R5}
Monodromies and punctures are determined not by curves, but by
Fuchsian equations. Therefore punctures, as attributes of functions
and Poincar\'e $\tau$-domains of their automorphisms, are not bound
to be branch places of a certain curve. They may be located at any
places on any curves,  including even non-branch (regular) places.
This is just the case of torus \eqref{TOR}, whereas the function
$z(\tau)$ itself and its equation \eqref{Qz} may uniformize many
other curves with branch places $z=E_k$, \eg, the curves of the form
\begin{equation}\label{Wab}
v^m=z^n\,(z-1)^p(z+1)^q(z-3)^r(z+3)^s\,.
\end{equation}
Here, the numbers $\{m,n,p,q,r,s\}$ are allowed to be any complex
numbers so that we construct the formal parametrizations of
\emph{non-algebraic} dependencies (genus is not finite) by
monodromies of finite topological genus! In this regard even
$\tau$-forms for non-algebraic  but single-valued solutions
$y(\tau)$ to equation \Psix\ itself provide a large number of such
examples. Apart from PH-solutions \eqref{algPic}, \eqref{simp} here
is a simplest one (new?): for arbitrary $s\in\mathbb{C}$ the
function $y=x^s$ solves \eqref{p6} when
$(\alpha,\beta,\gamma,\delta)=\big(0,0,s^2,(s-1)^2-\frac12 \big)$.
See also comments as to work by Guzzetti \cite{guz} in Sect.~6 of
work \cite{br2}. A direct check shows that the self-suggested
generalization $y=x^s(x-1)^r$ does not exist for
$s,r\notin\mathbb{Q}$.
%There is a very simple analog of this phenomenon in the case of
%classical uniformization by elliptic functions. The Weierstrass
%equation
%\begin{equation}\label{wpp}
%\wp'^2=4\,(\wp-e_1)(\wp-e_2)(\wp-e_3)
%\end{equation}
%is equivalent to equation
%\begin{equation}\label{wpp2}
%y^2=4\,(\wp-e_2)(\wp-e_3)\,,
%\end{equation}
%where
%$$
%y=\frac{\wp'}{\sqrt{\wp-e_1}}\,.
%$$
%On the other hand, function $\sqrt{\wp-e_1}$ is a single-valued one
%as it is expressed through Weierstrass' functions $\sigma_{\!1}^{}$
% and
%$\sigma$ \cite{we2}. Namely,
%$\sqrt{\wp-e_1}\sim\sigma_{\!1}^{}/\sigma$. Hence it follows that
%transcendental function $\wp$ parametrizes not only nonelementary
%elliptic curve \eqref{wpp} but elementary zero genus curve
%\eqref{wpp2}
%as well. The point $e_1$ is a `superfluous' singularity
% of function
%$\wp$ on this curve. Another form of the aforesaid is the well-known
%elementary ($g=0$) relation $\mathrm{sn}^2+\mathrm{cn}^2=1$, whereas
%Jacobi's sn-, cn-functions are presented in parametrizations of
%elliptic ($g=1$) curves. Relation between \eqref{wpp}, \eqref{wpp2}
%is %the same as between our modular uniformizations \eqref{Wab} and
%\eqref{TOR}.
\end{remark}

Thus,  Painlev\'e \emph{PH-curves generate new universal
uniformizations forming towers of new curves}. The function
\eqref{z}, for example, uniformizes all the curves of the form
\eqref{Wab} in the sense that analytic function $v(\tau)$ determined
by the relations \eqref{z}, \eqref{Wab} is a single-valued function
in the entire domain of its existence, that is $\Hp$. We can even
enlarge this class by the change $z\mapsto \sqrt{z}$. This is
possible because of Picard's function $y(\tau)$ is a perfect square
and, from the first formula in \eqref{tor}, one follows that $z$ is
also the perfect square:
$$
z=-y\,\frac{(z-1)^2}{4\,u^2}\qquad\hence\qquad
z(\tau)=\left\{\big(z(\tau)-1\big)\,
\frac{\vartheta_3\,\theta_2}{2\,\vartheta_4\,\theta_1}\right\}^2
$$
(non-obvious fact for $\theta$-constant expression \eqref{z}).
Therefore Hauptmodul $r=\sqrt{z\,}(\tau)$ satisfies the equation
(Lemma~\ref{L1})
$$
[r,\tau]=4\,r^2\,\bcQ(r^2)+\frac32\,\frac{1}{r^2}
$$
(exercise: describe its singularities). We shall return to  example
\eqref{Wab} further below since one of its  particular cases is
related to the widely known  Jacobi's curves
\begin{equation}
\label{hyper}
v^2=z\,(z-1)(z-a)(z-b)(z-a\,b)\,.
\end{equation}
Appearance of such  curves  is not an exception and we observe in
passing that the way of getting the hyperelliptic formulae from
Picard--Hitchin's towers is  simpler than  the direct search for
hyperelliptic modular dependencies, say  of genus $g=2$, among
subgroups of known congruence groups of level $N$. In all the known
cases the level turns out to be very large; see for example analysis
of group $\bo\Gamma_0(50)$ in work \cite{birch}. We consider yet
another example because it is connected with very classical objects.

\subsection{Schwarz hyperelliptic curve\label{shwarz}}

The structure of Picard's solutions is such that we may apply to any
of them the substitution $\{x=p^4,$ $y=-p^2\,q\}$. Consider again
the small values of $N$ and carry out this substitution in the
second of the curves \eqref{uy} ($N=3$). We then obtain the genus
$g=3$ curve
\begin{equation}\label{pq}
4\,(p^2+p^{\sm2})\,q=q^4-6\,q^2-3\,,\qquad
\left\{p=\frac{\vartheta_4(\tau)}{\vartheta_3(\tau)}\,,\quad q=
\frac{\theta_2^2\big(\frac16\big|\tau\big)}
{\theta_1^2\big(\frac16\big|\tau\big)}\right\}.
\end{equation}
The structure of \eqref{pq} tells us that it should be hyperelliptic
and a simple calculation shows that it is isomorphic to the classical
Schwarz curve
\begin{equation}\label{x8}
\bo{y}^2=\bo{x}^8+14\,\bo{x}^4+1\,.
\end{equation}
This curve   repeatedly appears throughout both volumes of Schwarz's
\emph{Gesammelte Abhandlungen} \cite{schwarz} in different
contexts\footnote{Including the very 1867 work \cite[{\bf I:}
p.~13]{schwarz} wherein his famous Schwarz derivative $\{s,u\}$ had
arisen; original notation was $\bo{\mathit{\Psi}}(s,u)$.} but none
of parametrizations was known until recently. Birational
transformation between \eqref{pq} and \eqref{x8} is as follows
\begin{equation}\label{bir}
\begin{alignedat}{4}
p&=\frac{\bo{x}^4-\bo{y}-1}{4\,\bo{x}}\,,
\qquad&\bo{x}&=\frac{p\,(q^2-1)}{2\,(p^2+q)}\,,\\
q&=\frac{\bo{x}^4-\bo{y}+1}{2\,\bo{x}^2}\,,
\qquad&
\bo{y}&= \frac{q^2+3}{q^2\,(q^2-1)}\,(4\,p^2\,q+q^2+3)\,.
\end{alignedat}
\end{equation}
One  easily obtains parametrizations and Fuchsian equation for
$\bo{x}(\tau)$. Curiously, in doing so we arrive at further
enlargement of the tower for curve \eqref{x8}. This is a frequently
encountered situation in universal uniformization. Since the objects
are very classical, we state  the results as a separate proposition.

\begin{proposition}\label{P3}
The function $\bo{x}=\bo{x}(\tau)$ defined by
$\vartheta,\theta$-ratio \eqref{pq}, \eqref{bir} has a zero/pole
divisor determined by the quotient
\begin{equation}\label{xuni}
\bo{x}=\frac{\ded^3(\tau)\,\theta_2\big(\frac13\big|\tau\big)}
{\theta_1^2\big(\frac16\big|\tau\big)\,
\theta_3^2\big(\frac16\big|\tau\big)}
\end{equation}
and uniformizes the tower of curves formed by Burnside's curve
$\bo{z}^n=\bo{x}^5-\bo{x}$ and Schwarz's curve \eqref{x8}\/$:$
\begin{equation}\label{x8x5}
\bo{z}^2=
(\bo{x}^8+14\,\bo{x}^4+1)(\bo{x}^5-
\bo{x})\,.
\end{equation}
\end{proposition}

\begin{pf}
Derivation of the prime-forms \eqref{xuni} is just simplification.
In order to present $\bo{x}(\tau)$ defined by \eqref{bir} in form of
such a ratio one needs to use the duplication formula
$\theta_2^4(z)-\theta_1^4(z)=\vartheta_2^3\cdot\theta_2(2z)$,
standard quadratic  $\theta$-identities
\cite{abramowitz,bateman,weber}, and second relation in
\eqref{difalg}.

The sequential algebraic changes of variables $x\mapsto p\mapsto
\bo{x}$ defined by \eqref{bir} transform (Lemma~\ref{L1}) equation
\eqref{k2} into the following Fuchsian equation:
\begin{equation}\label{Qx85}
\begin{aligned}
{}[\bo{x},\tau]&= -\frac12\,\frac
{\bo{x}^{24}-102\,\bo{x}^{20}+
1167\,\bo{x}^{16}+
1964\,\bo{x}^{12}+
1167\,\bo{x}^8-102\,\bo{x}^4+1}
{(\bo{x}^8+14\,\bo{x}^4+1)^2
(\bo{x}^5-\bo{x})^2}\\
&=-\frac38\,\sum_\beta\frac{1}{(\bo{x}-\beta)^2}-
\frac12\,\sum_\alpha\frac{1}{(\bo{x}-\alpha)^2} +
\frac{(\bo{x}^4+7)\,(5\,\bo{x}^4-1)\,\bo{x}^3}
{(\bo{x}^8+14\,\bo{x}^4+1)
(\bo{x}^5-\bo{x})}\,,
\end{aligned}
\end{equation}
where summations run over roots of polynomials
$\beta^8+14\,\beta^4+1$ and $\alpha^5-\alpha$. Bearing in mind
Remark~\ref{R5}, this provides uniformization of curves \eqref{x8}
and \eqref{x8x5}. The $\tau$-behavior of $\bo x(\tau)$ at all the
branch or `superfluous regular' places is determined by coefficients
$\frac12$ and $\frac38$ in \eqref{Qx85}.\hfill $\blacksquare$
\end{pf}

\begin{example}\label{E3'}
We may, independently of birational transformations \eqref{bir} but
with usage of Theorem~\ref{T5}, derive that function \eqref{xuni}
satisfies an equation with rational coefficients and this equation
is \eqref{Qx85}.
\end{example}

The polynomial expressions above, written in homogeneous form
$z_1^8+14\,z_1^4\,z_2^4+z_2^8$, $z_1z_2(z_1^4-z_2^4)$, are widely
known as Schwarz--Klein's ground forms representing  symmetry group
of octahedron \cite{schwarz,fricke,hitchin2}. We thus obtain that
correlation of two simplest PH-curves \eqref{uy} gives a
uniformizing $\tau$-representation for their product.

\begin{example}\label{E4}
This is a nontrivial exercise to find 13-generated monodromy
representation $\G_{\bo{x}}=\langle T_k\rangle$ for the
Schwarz--Burnside equation \eqref{Qx85} with use of
$\vartheta,\theta$-representation \eqref{xuni} and Theorem~\ref{T6}.
Since equation is solvable the basis of corresponding solutions
$\Psi_{1,2}(\bo{x})$ attached to  Hauptmodul \eqref{xuni} is also
computable.
\end{example}

\subsection{Some conclusions\label{7.3}}

Summarizing the preceding material, we may formulate the following
general recipe:
\begin{itemize}
\item

Every algebraic curve \eqref{F} may be thought of as an
algebraic substitution/trans\-for\-ma\-tion $x\mapsto y$ in any
Fuchsian equation $[x,\tau]=\bcQ(x,y)$ with correct accessory
parameters. Conversely, any rational/algebraic substitution
$x\mapsto z{:}\;x=\Phi(z)$ (or $\Phi(x,z)=0$) may be thought of
as an algebraic curve generating the new curve $\Xi(y,z)=0$.
Suppose the \emph{local} ramifications $z(x)$ are compatible
with the local behavior $x=E+a\,\tau^n+\cdots$ in the sense that
$z=E'+b\,\tau^m+\cdots$, where $m\in\mathbb{Z}$. Then the
\emph{global} monodromy $\G_z$ (and $\G_y$ of course) has finite
genus; accessory parameters in the proper Fuchsian equation
$[z,\tau]=\widetilde\bcQ(z,y;x)$ are also correct. The universal
uniformization is characterized by that it `absorbs' the
arbitrary ramification orders; its tower is thus always
infinite.
\end{itemize}
At this point, it is worth noting that the `trivial' (zero genus)
rational substitutions like $x=R(z)$ may lead to nontrivial curves
of high genera and, on the other hand,  nontrivial substitutions
$\Phi(x,z)=0$ of higher genera may preserve the trivial genus of
both the monodromies  $\G_x$ and $\G_z$, or $\G_y$; see \cite{br5}
for explicit examples. Insomuch as all these transitions constitute
just substitutions, all the arising Fuchsian equations are
integrable one through another. In particular, we have the following
conclusion.

\begin{theorem}\label{T8}
Fuchsian equations \eqref{Qxy} corresponding to any algebraic
Painlev\'e uniformizing functions $y(\tau)$ are pullbacks of
hypergeometric ${}_2F_1$-equations by rational or algebraic
$($computable\/$)$ functions. These ${}_2F_1$-integrabilities form
an equivalence relation; \eg, equivalence with equation defining the
series $_2F_1\!\big(\frac12,\frac12;1\big|z\big)$ or Legendrian
$\ellK$ and $\ellK'$.
\end{theorem}

\begin{pf}
The $x(\tau)$-function satisfies Fuchsian Eq.~\eqref{k2} with three
singularities. Therefore it is solved through the hypergeometric
$_2F_1\!\big(\frac12,\frac12;1\big|z\big)$-series inversions.
Algebraic dependency \eqref{F} and common $\tau$ entail that the
$_2F_1$-integrability of \eqref{k2} lifts to the same kind one for
\eqref{Qxy}; except algebraic dependencies in $z$-arguments of
$_2F_1(a,b;c|z)$-functions. Indeed, the $\Psi$-equation for
$y$-variable comes from substitution \eqref{F} into \eqref{fuchs}.
Since any two algebraic Painlev\'e solutions $y_1(\tau)$,
$y_2(\tau)$ are algebraically related through $x$-variable, their
$\Psi$-equations are solvable one through another. We get algebraic
pullbacks between different $_2F_1$-functions. Again, thanks to
common $x$, the pullbacks satisfy the symmetry, reflection, and
transitivity properties; and define thereby an equivalence relation.
Representations~\eqref{kk'} show that all the integrabilities may be
represented in terms of Legendrian functions $\ellK$ and
$\ellK'$.\hfill $\blacksquare$
\end{pf}

It should be noted here that multi-valuedness (coming from
algebraicity) of $z$-argu\-ments above is not essential in $_2F_1$
because the hypergeometric function itself is multi-valued and
defined not uniquely. There are numerous transformations between
different $_2F_1$-objects \cite{abramowitz,bateman,goursat}. Also,
it does not matter which kind of multi-valuedness we meet:
algebraic, $_2F_1$-transcendental, or both of them. One finally
needs only the single-valued inversions of \eqref{ratio}.

Thus, \emph{uniformization is always split into uniformizations on
the equivalence classes} and the theory for Eq.~\eqref{p6} forms its
proper closed class: the \Psix-class. Corresponding monodromies $\G$
are aligned in towers. In the next two sections we shall exhibit all
the kind integrability equivalencies: rational/algebraic, between
PH' and \hbox{non-PH'} cases.

\section{Relations to Painlev\'e non-PH-curves\label{nonPH}}

Our constructions considered up until this point were dealing with
PH-curves being Painlev\'e curves or their consequences being no
solutions to \Psix. We can, however, get interesting information
involving Painlev\'e curves being no non-PH-ones.

Consider  a genus zero Painlev\'e non-PH-curve $F(x,y)=0$. It has
some rational parametrization $x=R_1(\T)$, $y=R_2(\T)$. On the other
hand, by virtue of universality of $x(\tau)$, it has parametrization
$y=y(\tau)$ through the `universal' $\tau$, wherein $y(\tau)$ is as
yet unknown. Rational parametrization through the $\T$ tells us that
$\T$ itself is a rational function of $x,y$, that is $\T=R(x,y)$.
Hence $\T$ becomes a univalent function of $\tau$ and therefore must
satisfy a Fuchsian \ode\ with the rational $\bcQ(\T)$-function:
$$
\T=R\big(x(\tau),y(\tau) \big)\qquad\hence\qquad
[\T,\tau]=\bcQ(\T)\,.
$$
In general, this way leads to \emph{new} universal Hauptmoduln
$\T(\tau)$ with monodromies known to be Fuchsian. We may further
take one more curve $\tilde F(x,\tilde y)=0$ (\eg, PH-curve) with a
known Hauptmodul $\tilde y(\tau)$. Rational uniformizer $\tilde \T$
for this second curve and the old one $\T$ do certainly lie on a
certain curve $\Xi(\T,\tilde \T)=0$; it may have, however, quite
high genus. If, on the other hand, $\tilde \T$ is a linear
fractional function of $\T$ then $y$ becomes a rational function on
the second curve. This is just we want because we may apply all the
preceding machinery since we obtain the rational function on the
(known) PH-curve but the function itself is a field generator for
the new (non-PH's) curve. If genus of $\Xi$ is nontrivial we obtain
new nontrivial $_2F_1$-integrable Fuchsian equations with correct
accessory parameters; we recall Sect.~\ref{7.3}. Let us illustrate
the aforesaid. The examples that follow can be enlarged from a large
collection of solutions listed in works \cite{boalch},
\cite{boalch2}, and \cite{lisovyy}.

\begin{example}\label{E7}
Consider a non-PH-solution obtained by P.~Boalch in work
\cite{boalch}:
\begin{equation}\label{81}
y=\frac{7\,\T^2+22\,\T+7}{8\,(\T^2+\T+1)\,(\T+2)\,\T}\,,\qquad
x=\frac{2\,\T+1}{(\T+2)\,\T^3}\,.
\end{equation}
Correlate it with Picard's solution \eqref{uy} (changing there
$y\mapsto \tilde y$) which is parametrized as follows
\begin{equation}\label{82}
\tilde y=-3\,\frac{(\tilde \T-3)\,(\tilde \T+1)}{(\tilde
\T+3)^2}\,,\quad x=\frac{(\tilde \T+1)(\tilde
\T-3)^3}{(\tilde \T-1)\,(\tilde \T+3)^3} \qquad\hence\quad\tilde\T=
\frac{\tilde y^3-3\,x\,\tilde
y+x\,(x+1)}{x\,(x-1)}\,.
\end{equation}
Equating the $x$-parts of \eqref{81} and \eqref{82} to each other,
we get $\tilde\T=\frac{3\Smaller[3]{\mathrm{T}}-1}
{\Smaller[3]{\mathrm{T}}-1}$. Therefore rational function
$\T=R(x,y)$ on Boalch's curve \eqref{81} (expression is too large to
display here) satisfies the following Fuchsian equation
\begin{equation}\label{83}
[\T,\tau]=
\frac{-2\,(\T^2+\T+1)^4}{(\T^2-1)^2\,(2\,\T+1)^2\,(\T+2)^2\,
\T^2}\,.
\end{equation}
The simplest way to obtain this equation is to apply Lemma~\ref{L1}
to \eqref{k2} with the second formula in \eqref{81}. Equation
\eqref{83} defines universal uniformizing Hauptmodul for `parabolic'
points $\T=\mbig\{0,\pm1,-2,-\frac12,\infty \mbig\}$, but this
Hauptmodul is in fact not new because transformation $\T\mapsto z$
of the form
$$
z=3\,\frac{\T+1}{\T-1}
$$
turns \eqref{83} into \eqref{Qz} so that we come back to the known
consequence of Picard's curve \eqref{uy}. We obtain, however,
Boalch's $y$ as a rational function on Picard's curve \eqref{82}:
$$
y=\frac{1}{16}\,\frac{15\,\tilde y^3-(14\,\tilde y^2+3\,\tilde
y-18)\,x}{(\tilde y^2-3\,\tilde y+3)\,\tilde y}\,.
$$
This gives uniformizing $y(\tau)$-representation \eqref{subswp}
which is not obvious a priori. Indeed, uniformization of algebraic
non-PH-solutions is as yet unknown because it does not follow
automatically from the PH-series.
\end{example}

In other words, the `hyperelliptic' Hauptmodul $z(\tau)$ defined by
\eqref{z}, \eqref{Qz}, or by Proposition~\ref{P5} is not merely a
nice example but exhibits some general recipe.

\begin{itemize}
\item Hauptmoduln coming from the PH-series link the zero genus
    Painlev\'e curves related to each other through Okamoto's
    transformations that, however, \emph{\/fall outside the
    scope of pure Picard--Hitchin's transformations \eqref{ok}}.
\end{itemize}
This construction can be continued with generating the new curves.
Here is a less simple example leading to a non-rational curve.

\begin{example}\label{E8}
Consider parametrization of the three-branch tetrahedral Hitchin's
solution corrected by Boalch as it has been written in \cite[formula
(10)]{boalch}:
$$
y=\frac{(\T-1)\,(\T+2)}{(\T+1)\,\T}\,,\qquad
x=\frac{(\T-1)^2\,(\T+2)}{(\T+1)^2\,(\T-2)}\,.
$$
Parameters for this solution are as follows
$$
\alpha=\frac19\,,\quad
\beta=4\,c^2\,,\quad \gamma=c^2\,,\quad \delta=c^2-\frac12\,,
$$
where $c$ is free. Formula for $x(\T)$ and Lemma~\ref{L1} imply
$$
[\T,\tau]=-\frac12\,\frac{\T^6+6\,\T^4-15\,\T^2+44}
{(\T^2-1)^2(\T^2-4)^2\,\T^2}\,.
$$
This Hauptmodul is of course universal but corresponds to the set of
five  parabolic singularities $\T=\{\pm1,\pm2,\infty \}$ rather than
six ones. The curve $\Xi(\T,\tilde \T)=0$ binding Picard's
$\tilde\T$ \eqref{82} and  $\T$ turns out to be an elliptic curve
and we easily find
$$
\Xi\colon\qquad \tilde\T^4+4\,(\T^3-3\,\T)\cdot\tilde\T^3
+18\, \tilde\T^2-27=0\,.
$$
Corresponding Picard's and Boalch's solutions  also lie on the
nontrivial elliptic curve
$$
(4\,\tilde y-3)(\tilde y^4+6\,\tilde y-3)\,y^6 - 12\,\tilde y\,
(4\,\tilde y-3)(\tilde y^3+3\,\tilde y-2)\,y^5+\cdots- 12\,\tilde
y^5\,(\tilde y^3+3\,\tilde y-2)\,y+4\,\tilde y^8=0\,,
$$
where we reduced the formula and designated by dots the polynomials
in descending powers of $y$. As above one can find Boalch's $y$ and
Picard's $\tilde y$ as rational functions on Picard's and Boalch's
curves respectively. These two curves are  of course isomorphic and
easy computation shows that they are birationally equivalent to the
simple Weierstrass canonical form $v^2=4\,u^3-48\,u+80$. Although
this torus and torus \eqref{tor-yu} arise  form one PH-curve, they
have different $J$-invariants: $-\frac{16}{9}$ and
$\frac{2197}{972}$ respectively.

The symmetry $\T\mapsto -\T$ prompts us to consider the Hauptmodul
$\bo{\T}=\T^2$ and we readily obtain its Fuchsian equation:
$$
[\bo{\T},\tau]=-\frac38\,\frac{1}{\bo{\T}^2}
-\frac12\,\frac{1}{(\bo{\T}-1)^2}
-\frac12\,\frac{1}{(\bo{\T}-4)^2}+ \frac18\,
\frac{7\,\bo{\T}-11}
{(\bo{\T}-1)\,(\bo{\T}-4)\,\bo{\T}}
$$
which defines a thrice  punctured sphere at $\bo{\T}=\{1,4,\infty
\}$ with an additional conical singularity of the second order at
$\bo{\T}=0$. Such equations can also be applied to uniformization
(mixed uniformization in terminology of \cite{br}). Solution to
inversion problem for  $\T$ and $\bo{\T}$ is somewhat nonstandard
and we shall present it elsewhere.
\end{example}

\begin{example}\label{E9}
Consider any PH-curves for $N=\{3,4\}$. Then introduction of the
parameter $\T$ as above leads to a  nice Fuchsian \ode
$$
[\T,\tau]=-\frac12\,\sum\limits_\alpha\frac{1}{(\T-\alpha)^2}+
\frac{4\,\T^6-14\,\T^4+2}{\T^9-6\,\T^7+6\,\T^3-\T}\,,
$$
where $\alpha^9-6\,\alpha^7+6\,\alpha^3-\alpha=0$. Solutions to this
\ode\ are expressible in terms of  $_2F_1$-series. Universal
Hauptmodul and its differentials can also be derived. One easily
computes that introduction of the second parameter $\tilde\T$ leads
to a (uniformizable) curve $\Xi(\T,\tilde\T)=0$. It has genus $g=7$
and some elements of its tower can also be analyzed. Notice, genera
of the initial PH-curves are still equal to zero (see
table~\eqref{genera}).
\end{example}

\section{Inversion problems and related topics\label{Fuchs}}

Uniformization of  PH-curves leads to many remarkable consequences
and there is no escape from the mentioning appearance of the two
famous Fuchsian equations in the `algebraic \Psix-theory'.

\subsection{Exceptional cases of Heun's equations $($Chudnovsky
equations\/$)$\label{Heun}}

The symmetry $z\mapsto -z$ of equation \eqref{Qz} suggests to make
the substitution $z^2=s$ and to expect some simplification of this
equation. Doing this, we obtain that function $s(\tau)$ satisfies
the \ode
$$
[s,\tau]=-\frac12\,{\frac{s^4- 12\,s^3+102\, s^2-108\,s+
81}{s^2\,(s-1)^2\,(s-9)^2}}
$$
which, in the language of the `linear Fuchs theory', corresponds to
\begin{equation}\label{heun}
\Psi''=-\frac14\!\left(\,\frac{1}{s^2}+ \frac{1}{(s-1)^2}+
\frac{1}{(s-9)^2}- \frac{2\,s-8}{s\,(s-1) (s-9)}\right)\!\Psi\,.
\end{equation}
It is a Heun equation with singularities at $s=\{0,1,9,\infty\}$
defining  the monodromy of a 4-punctured sphere since exponent
differences at all the singularities are equal to zero. It is
notable that this equation is not treated by any classical
algorithmic methods in theories of integration of linear \odes\ (the
Picard--Vessiot theory \cite{singer}). This fact is not strange.
Presently only four exceptional  cases of integrable Fuchsian
equations with such a type of parabolic singularities are known.
They were claimed by D.~Chudnovsky \& G.~Chudnovsky in work
\cite{chud2} and read in their notation as follows \cite[p.~185;
correcting a typo for case \eqref{III}]{chud2}
\begin{alignat}{2}
&x\,(x^2-1)\,y''+(3x^2-1)\,y'+x\,y=0\,,\tag{I}\label{I}\\
&x\,(x^2+3x+3)\,y''+(3x^2+6x+3)\,y'+(x+1)\,y=0\,,\tag{II}\label{II}\\
&x\,(x-1)(x+8)\,y''+(3x^2+14x-8)\,y'+(x+2)\,y=0\,,
\tag{III}\label{III}\\
&x\,(x^2+11x-1)\,y''+(3x^2+22x-1)\,y'+(x+3)\,y=0\,.\tag{IV}\label{IV}
\end{alignat}
The case \eqref{I} is equivalent to equation \eqref{lemn}; we have
detailed it in {\sf Example~\ref{E2}}. The case \eqref{III} becomes
equation \eqref{heun} after the change
\begin{equation}\label{sx}
x=1-\frac9s\,.
\end{equation}

Explicit integrability of equation \eqref{heun}  follows immediately
from Sect.~\ref{Hyper}. Since $s$ is algebraically related to
$u=\sqrt{x}=k'(\tau)=k\mbig(\frac{-1}{\tau}\mbig)$, we obtain from
\eqref{uy} and \eqref{tor}:
\begin{equation}\label{xs}
64\,(2\,x-1)^2\,s=(s^2-6\,s-3)^2\,.
\end{equation}
Therefore this algebraic change $s\mapsto x$, supplemented with the
linear change $\Psi\mapsto\psi$:
\begin{equation}\label{tmp}
\Psi=\sqrt{s_x}\cdot\psi(x)\qquad\hence\qquad
\Psi=\frac{\sqrt[\uproot{1}\leftroot{1}4]{s^3\,}}{s-1}
\cdot\psi(x)\,,
\end{equation}
reduces  \eqref{heun} to  a hypergeometric equation in normal form
\eqref{k2}:
$\psi^{}_{\mathit{xx}}=-\frac14\,\frac{x^2-x+1}{x^2(x-1)^2}\,\psi$.
We get solution:
\begin{equation}\label{sK}
\Psi=\sqrt[\uproot{1}4]{s\,(s-1)^2(s-9)^2}\cdot
\mbig\{A\,\ellK(\sqrt{x})+B\,\ellK'(\sqrt{x})\mbig\}\,,
\end{equation} where
$$
x=\frac{1}{16}\,\big(3\!{\ts\sqrt[-2]{s\,}}+1\big)
\big(\ts\sqrt{s\,}-1\big)^3\,.
$$
Recall that symbols $\ellK$ and $\ellK'$ represent hypergeometric
functions \eqref{kk'} and all their analytic continuations.

It is clear that further examples are rapidly multiplied and, what
is more important, for all equations  we have always had solutions
to the inversion problems since these solutions are \emph{hidden}
forms of the one basic relation \eqref{tau} (Theorem~\ref{T8}). For
example, for Heun's equation \eqref{heun} inversion of the ratio of
its two solutions
\begin{equation}\label{invs}
\frac{\Psi_1}{\Psi_2}=\frac{\ellK\mbig[5](\mfrac14\sqrt{(s-6)
\sqrt{s\,}-
3\!\sqrt[\leftroot{-1}\uproot{1}-2]{s\,}+8\,}\,\mbig[5])} {{\ds
\ellK'}\mbig[5](\mfrac14\sqrt{(s-6)\sqrt{s\,}-
3\!\sqrt[\leftroot{-1}\uproot{1}-2]{s\,}+8\,}\,\mbig[5])}=
\frac{a\,\tau+b}{c\,\tau+d}\,,
\end{equation}
that is function $s=s\mbig[3](\frac{a\,\tau+b}{c\,\tau+d}\mbig[3])$,
is a square of the $\theta$-constant expression \eqref{z} with an
appropriate choice of constants $(a,b,c,d)$. In order to determine
them it is sufficient to consider any three points $\tau$
corresponding to any tori with complex multiplication since such
tori have exact (algebraic over $\mathbb{Q}$ \cite{weber}) values of
$\ellK$. The function $s(\tau)$ thus constitutes one further
universal uniformizing function for a new set of points. Any
algebraic function of $s(\tau)$ with arbitrary ramifications at
points $s=\{0,1,9,\infty\}$ is a single-valued function of $\tau$.
The complete solution to inversion problem \eqref{invs} and
$\tau$-representation for the $\Psi$-function will be given in
Sect.~\ref{further}.

It should be noted here that solution form depends on which group
(variable) is chosen\footnote{The author is indebted to M.~van~Hoeij
for (perhaps) most elegant rational $_2F_1$-pullback in this case:
$$\Psi(s)=\!\!\!{\sqrt[\sm6]{s-1}}\:\sqrt{s\,(s-9)}\:
{}_2F_1\bigg(\frac13,\frac13;1\bigg|{-}
\frac{s\,(s-9)^2}{27\,(s-1)^2} \bigg).$$} and these different
choices entail some byproducts. For instance,  involving the group
$\bo\Gamma(1)$, and therefore Klein' invariant $J$, we know that
$$
J=\frac{4}{27}\frac{(x^2-x+1)^3}{x^2\,(x-1)^2}\,.
$$
Correlating this expression with  formula \eqref{xs}, we get one
more integrating change\ and  nonstandard representation for Klein's
$J(\tau)$ through the $\vartheta,\theta$-constants \eqref{z}.

\begin{proposition}\label{P4}
Klein's invariant $J(\tau)$ has the following $\theta$-constant
representation
\begin{equation}\label{J}
J=\frac{1}{12^3}\,
\frac{(s+3)^3\mbig[1]((s-5)^3+2^7\mbig[1])^3}{s\,(s-1)^6(s-9)^2}\,,
\end{equation}
where $s=s(\tau)$ is a square of expression \eqref{z} $($see also
formula \eqref{Ztau} below\/$)$.
\end{proposition}

Substitution \eqref{J} solves \eqref{heun} in terms of
hypergeometric functions $_2F_1\big(\frac{1}{12},\frac{1}{12};
\frac23\big|\,J\big)$. Recall, this is also true for arbitrary
Painlev\'e curves at all. According to Theorem~\ref{T8}, any
Fuchsian equation arising from Painlev\'e substitution
\eqref{substheta} is integrable in terms of $\ellK,\ellK'(\sqrt{x})$
or $_2F_1(J)$-series. These series are directly related to the
classical representations for $J(\tau)$ through $\vartheta(\tau)$'s
or Dedekind's function $\ded(\tau)$. It is known that the
$\ded(\tau)$ itself is a hidden form of the $\theta$-constant
$$
\ded(\tau)=-\ri \re^{\frac{\pi}{3}\ri\tau}_{\mathstrut}\,
\theta_1(\tau|3\tau)
$$
and is a particular case of the new $\theta$-constant class
described in sects.~\ref{diffcalc} and \ref{groups}. There are of
course many other analogs of \eqref{J} and therefore PH-curves and
their consequences provide a further and rich development of
theories to the classical $\vartheta$- and $\ded$-constants. As for
solutions to the $\Psi$, Picard--Hitchin's hierarchy leads to a
massive generalization of numerous quadratic and cubic
transformations of $_2F_1$-series listed in dissertation by Goursat
\cite{goursat}. In addition to this, we have also representations of
groups, $\theta$-Hauptmoduln, their differential calculus, etc. The
aforesaid and examples are the illustrations to what we said in the
end of Sect.~\ref{universal}: substitutions bind integrable
equations between themselves. For exhaustive theory to algebraic
equivalence of the Chudnovsky list see \cite{br5}.

\subsection{Ap\'ery's differential equations\label{apery}}

Roger Ap\'ery, in his celebrated proof \cite{apery,poorten} of the
irrationality of Riemann's
$\zeta(3)=1^{\sm3}+2^{\sm3}+3^{\sm3}+\cdots$, used the recursion
$$
n^3\,C_n=(34\,n^3-51\,n^2+27\,n-5)\,C_{n-1}-(n-1)^3\,C_{n-2}
$$
and pointed out that  it corresponds to  the linear 3rd order \ode
\begin{equation}\label{A3}
r^2(r^2-34\,r+1)\,\psi'''+r\,(6\,r^2-153\,r+3)\,\psi''+
(7\,r^2-112\,r+1)\,\psi'+(r-5)\,\psi=0
\end{equation}
with the help of standard correlation between $C_n$-recursions and
solutions to  linear \odes: $\psi=\sum C_n r^n$. He also observed
that this \ode\ was a second symmetric power of the following  \ode\
of Fuchsian class (see also \cite{dwork}):
\begin{equation}\label{A2}
r\,(r^2-34\,r+1)\,\varphi''+(2\,r^2-51\,r+1)\,\varphi'+
\frac{1}{4}\,(r-10)\,\varphi=0\,.
\end{equation}

Uniqueness of Ap\'ery's equation \eqref{A3} suggests to seek for its
relatives among some `good` equations, \ie,  equations with known
monodromies, $_2F_1$-integrals, etc. By virtue of pointed out
relation between \eqref{A3} and \eqref{A2} we may restrict our
consideration to the 2nd order equation \eqref{A2} since solution of
\eqref{A3} is given by the formula
$\psi=a\,\varphi_1^2+b\,\varphi_1^{}\varphi_2^{}+c\,\varphi_2^2$.
Without entering into details of irrationality, integrality of
Ap\'ery's sequences, etc,  we note that relation of these sequences
to modular forms for certain subgroups of $\bo\Gamma(1)$ was
established in the beautiful paper of Beukers \cite{beukers}. We
shall give an independent motivation/derivation and explanations in
the context of Picard--Hitchin's uniformization.

Let us transform \eqref{A2} into the normal (and unique) Klein's
form $\varphi''=\frac12\,\bcQ(r)\,\varphi$ by the linear change
$\varphi\mapsto \sqrt[4]{r^4-34\,r^3+r^2}\:\varphi$. We get Heun's
equation
\begin{equation}\label{A2Klein}
\varphi''=-\bigg\{\frac12\,\frac{1}{r^2}+\frac38\,
\sum_\alpha\frac{1}{(r-\alpha)^2}
-\frac34\,\frac{r-16}{(r^2-34\,r+1)\, r}\bigg\}\varphi\,,
\end{equation}
where $\alpha$'s are roots of equation $\alpha^2-34\,\alpha+1=0$.
Exponent differences $\delta$ for this equation are $\delta=0$ at
points $r=\{0,\infty\}$ and $\delta=\frac12$ at $r=\alpha$'s.

Suppose we do nothing know about solutions to Eq.~\eqref{A2Klein}.
Motivated by a desire to reduce it to some integrable form, we
should try a transformation $r\mapsto s$ that must be of
algebraic/rational type lest the Fuchsian class be escaped. The
one-to-one linear fractional transformation will nothing yield and
we try the rational one at first: $r=R(s)$. Obviously, images of two
parabolic points $r=\{0,\infty\}$ will remain parabolic ones in
$s$-equation under such a transformation since at these points we
have
\begin{equation}\label{ln}
\frac{\varphi_1^{}}{\varphi_2^{}}\sim\ln r+\cdots=\ln R(s)+\cdots;
\end{equation}
so at least two parabolic singularities are not removable.
Considering other singularities, we conclude that if the
transformation $r\mapsto s$ were  regular at $r=\alpha$, it would
cause the total number of $s$-singularities to increase. Indeed,
$$
\frac{\varphi_1^{}}{\varphi_2^{}}\sim\sqrt{r-\alpha}+\cdots=
\sqrt{R(s)-\alpha}+\cdots
$$
and analytic function $R(s)$ has greater than two $\alpha$-points
for each $\alpha$. New Fuchsian $s$-equation would be one with
greater number of singularities instead of 4 points in initial
Eq.~\eqref{A2Klein}; it is undesirable. Therefore we want each
$r=\alpha$ gets mapped into one $s$-point and the points themselves
should be points of regularity in $s$-equation.

All the $s$-singularities arise only from \eqref{ln}; they will
inevitably be again parabolic. The resulting $s$-equation will thus
determine a punctured sphere. So, the presence of two conical points
$r=\alpha$ in \eqref{A2Klein} tells us that
$(r\,{\mapsto}\,s)$-transformation  must be a quadratic one of the
form
\begin{equation}\label{rs0}
r=\frac{(a\,s-b)(s+c)}{(s-d)(s-e)}
\end{equation}
and its $s$-discriminant $\Delta$ should read as follows:
\begin{equation}\label{discrim}
\Delta=r^2-34\,r+1\,.
\end{equation}

Three $s$-images of singularities  can be freely appointed so that,
without loss of generality, we may assign $(r=\infty)\mapsto
(s=\{0,1\})$ and $(r=0)\mapsto (s=\infty)$ in \eqref{rs0}. Other
$s$-images may not be `stirrable' as we have the discriminant
condition \eqref{discrim}:
$$
\left\{r=\frac{a\,s-b}{s\,(s-1)},\quad\Delta=r^2-34\,r+1\right\}\,.
$$
Hence $\{a^2=1,\;2\,b-a=17\}$ and we thus get
\begin{equation}\label{rs}
r=\frac{s-9}{s\,(s-1)}\qquad
\text{or}\qquad r=\frac{s+8}{s\,(1-s)}\,.
\end{equation}
The $s$-variable is chosen up to a linear fractional transformation
and we easily find that these two solutions are in fact the only one
since the changes
$$
s\mapsto1-s\,,\qquad s\mapsto\frac{9\,s}{s+8}
$$
turn one solution in \eqref{rs} into another.

It is a remarkable fact that the resulting (unique) solution leads
exactly to \eqref{heun}; Ap\'ery's equation turns out to be a hidden
consequence of the two simplest Dubrovin--Hitchin solutions
\eqref{uy} coming from an infinite series of Picard's ones. Indeed,
the change $\varphi(r)=\sqrt{r_s}\:\Psi(s)$ and first of equalities
\eqref{rs} substituted in  \eqref{A2Klein} cause this equation to
become Eq.~\eqref{heun}. The second equality produces  an equivalent
equation. This and other equivalents are obtained from \eqref{heun}
by the six linear fractional transformations $s\mapsto
\frac{a\,s+b}{c\,s+d}$ permuting points $s=\{0,1,\infty\}$ between
themselves and therefore the fourth singularity can be freely
appointed to one of the six values
$s=\big\{9,\frac19,\frac89,\frac98,-8,-\frac18 \big\}$. Formulae
\eqref{rs} lead to solution of inversion problem for Ap\'ery's
equation \eqref{A2} in terms of those for $s$-equation \eqref{heun}.

\begin{remark}
When passing to a punctured sphere above, we had actually motivated
there an important transition from non-free monodromy $\G_r$ to a
free one $\G_s$ and obtained thereby yet another explanation as to
why universal uniformization by free groups leads to large `towers
of solvable curves'. As a rough guide one may think that \emph{free}
monodromy, being a rather large group, integrates many other
equations but not only its proper one. The \emph{larger} groups are
\emph{easier} described. For example group $\G_s$ is larger than
$\G_r$ because it is completely free of defining relations, whereas
$\G_r$ has the two ones $\mathfrak{a}^2=\mathfrak{b}^2=1$. Number of
their generators is the same.
\end{remark}

\subsection{Hauptmoduln and the $\Psi$ $($examples\/$)$. Is the
calculus closed?\label{further}}

Having a structure description  of Heun--Ap\'ery's equations
\eqref{heun}, \eqref{A2}--\eqref{A2Klein} we can now involve
simplifications with use of transformations of
$\vartheta,\theta$-constants \cite{weber,knopp} and sum up the
preceding stuff. We shall do this without entering into details of
calculations. For example, it becomes immediate to complete formula
uniformization for equation \eqref{heun} and related equations
\eqref{Qz}, \eqref{hyper}.

\begin{theorem}\label{T9}
The  square of Hauptmodul \eqref{z}
\begin{equation}\label{Ztau}
s(\tau)=9\,\frac{\vartheta_3^4(3\,\tau)}{\vartheta_3^4(\tau)}=
\frac{\theta_1^4\big(\frac13\big|\tau\big)\,
\theta_3^4\big(\frac13\big|\tau\big)}
{\theta_2^4\big(\frac13\big|\tau\big)\,
\theta_4^4\big(\frac13\big|\tau\big)}
\end{equation}
solves the inversion problem \eqref{invs} for Heun's equation
\eqref{heun}. The three group generators
\begin{equation}\label{star}
\G_s=\bigg\langle
T_{(0)}=\mbig[7](\,\begin{matrix}4&\!\!-1\\9&\!\!-2\end{matrix}\,
\mbig[7])_{\tau=\frac13},\quad
T_{(\infty)}=\mbig[7](\,
\begin{matrix}4&\!\!-3\\3&\!\!-2\end{matrix}\,\mbig[7])_{\tau=1},
\quad T_{(9)}=\mbig[7](\,\begin{matrix}1&\!\!-2\\0&
\!\!\phantom{-}1\end{matrix}\,
\mbig[7])_{\tau=\infty}\bigg\rangle
\end{equation}
determine the monodromy representation for  equation \eqref{heun}
and automorphism of \eqref{Ztau}. Here, subscripts indicate the fix
$\tau$-points of \ $T$'s and values of $s$-singularities at them.
The fourth point $s(0)=1$ corresponds to the cycle of cusps
$$
\tau=0\;\;\overset{T_{(0)}}{{\scalebox{1.7}[1]{$\mapsto$}}}\;\;
\frac12
\;\;\overset{T_{(\infty)}}{{\scalebox{1.7}[1]{$\mapsto$}}}\;\;2\;\;
\overset{T_{(9)}}{{\scalebox{1.7}[1]{$\mapsto$}}}\;\;0\,,
$$
so that\/ $T_{(9)}\,T_{(\infty)}\,T_{(0)}=\mbig[3](
\begin{smallmatrix}\phantom{-}1&0\\-6&1\end{smallmatrix}\mbig[3])=
T_{(1)}$.
\end{theorem}

Some consequences suggest themselves. From \eqref{xs}--\eqref{tmp}
we obtain (again, after $\vartheta,\theta$-simplifications) an
explicit $\tau$-representation of the function
$\Psi(s)=\tilde\Psi(\tau)$ since one may always normalize
$$
\tilde\Psi_2(\tau)=\sqrt{\dot s(\tau)\,}\,,
\qquad\tilde\Psi_1(\tau)=\tau\,\tilde\Psi_2(\tau)\,.
$$

\begin{corollary}\label{C4}
The Heun equation  \eqref{heun} has the following uniformizing
$\tau$-representation for its solutions $\Psi_{1,2}(s)$ attached to
Hauptmodul \eqref{Ztau}\/$:$
\begin{equation}\label{Psi-s}
\tilde\Psi_2(\tau)=\frac{\vartheta_3^3(3\,\tau)}{\vartheta_3(\tau)}\,
\frac{\vartheta_2^4\mbig(\frac{\tau+1}{2}
\mbig)}{9\,\vartheta_3^4(3\,\tau)-\vartheta_3^4(\tau)}\,,\qquad
\tilde\Psi_1(\tau)=\tau\,\tilde\Psi_2(\tau)\,.
\end{equation}
\end{corollary}

The relation \eqref{Psi-s} is an analog of the well-known
$\bo\Gamma(2)$-Jacobi $\tau$-representation for the elliptic
integral $\ellK(k)=\frac{\pi}{2}\,\vartheta_3^2(\tau)$. Notice
incidentally that this identity and commensurability of all the
groups in question in $\bo\Gamma(1)$ imply that all such identities
are merely hidden subgroup \mbox{$\tau$-representations} of the
general and well-known transformation law
$$
\ded^2(\tau)
\quad\overset{\widehat{\bo\Gamma(1)}(\tau)}
{{\scalebox{3.2}[1]{$\mapsto$}}} \quad
\Aleph^{\frac23}_{\mathstrut}\!\cdot(c\,\tau+d)\,\,\ded^2(\tau)\,.
$$
It is in turn a hidden form of solution to the known linear Fuchsian
equation with $\bo\Gamma(1)$-monodromy
$$
J(J-1)\,\psi''+\frac16(7J-4)\,\psi'+\frac{1}{144}\,\psi=0
$$
and therefore is a corollary of the obvious fact
$$
\psi
\quad\overset{\G_{\!\!\!\!\sss J}^{}\,=\,\bo\Gamma(1)}
{{\scalebox{3.5}[1]{$\mapsto$}}} \quad (c\,\tau+d)\,\psi\,.
$$
It follows  immediately that $\psi=\ded^2(\tau)$  (check this
directly).

\begin{remark}[\sf Example]\label{R7}
It is well known that differential calculus associated with
$\Psi$-functions for Legendre's modulus $k^2(\tau)$ contains the
closed set of functions $\{\ellK$, $\ellK'$, $\ellE$, $\ellE'\}$;
the derivatives $\frac{d \ellK}{dk}$, $\frac{d\ellE}{dk}$, \...\ are
functions of $\ellK$, $\ellE$'s themselves \cite{bateman},
\cite[{\bf IV}:~p.~157]{tannery}:
\begin{equation}\label{KE1}
\begin{aligned}
\frac{d\ellK}{dk}&=-\frac{\ellK}{k}-\frac{\ellE}{(k^2-1)\,k}\,,
&\qquad\quad
\frac{d\ellK'}{dk}&=\frac{k\,\ellK'}{1-k^2}+
\frac{\ellE'}{(k^2-1)\,k}\,,
\\
\frac{d\ellE}{dk}&=-\frac{\ellK}{k}+\frac{\ellE}{k}\,,&\qquad\quad
\frac{d\ellE'}{dk}&=\frac{k\,\ellK'}{1-k^2}+
\frac{k\,\ellE'}{k^2-1}\,.
\end{aligned}
\end{equation}
As a references source, we present here also the modular
representations for $\{\ellE,\ellE'\}$ since no complete set of such
formulae seems to have been tabulated in the standard texts (see
also comments and discussion as to classical elliptic integrals in
\cite[Sect.~7.3--4]{maier}). Formulae for $\ellK$ and $\ellK'$ are
known
\begin{equation}\label{KE2}
\ellK(k)=\frac{\pi}{2}\,\vartheta_3^2(\tau)\,,\qquad
\ellK'(k)=\frac{\pi}{2\ri}\,\tau\,\vartheta_3^2(\tau)\,,\qquad
k=\frac{\vartheta_2^2(\tau)}{\vartheta_3^2(\tau)}
\end{equation}
and formulae for $\ellE$, $\ellE'$ read
\begin{equation}\label{KE3}
\begin{aligned}
\ellE(k)&=\frac{2}{\pi}\frac{1}{\vartheta_3^2(\tau)}\!
\bigg\{\phantom{\tau}\,\eta(\tau)+
\frac{\pi^2}{12}\mbig{[}\vartheta_3^4(\tau)+\vartheta_4^4(\tau)\mbig]
\bigg\},\\
\ellE'(k)&=\frac{2}{\pi}\frac{\ri}{\vartheta_3^2(\tau)}\!
\bigg\{\tau\,\eta(\tau)-
\frac{\pi^2}{12}\mbig{[}\vartheta_2^4(\tau)+
\vartheta_3^4(\tau)\mbig]\,\tau
-\frac{\pi\ri}{2} \bigg\}.
\end{aligned}
\end{equation}
Lemma~\ref{L2} closes all the differential computations related to
the classical objects $\{\ellK$, $\ellK'$, $\ellE$, $\ellE'\}$ in
both $k$- and $\tau$-representations. In a reverse direction the
formulae read \cite[18.9.13]{abramowitz}
\begin{equation}\label{KE4}
\tau=\ri \frac{\ellK'(k)}{\ellK(k)}\,,\qquad\eta(\tau)=\ellK(k)
\,\ellE(k)+\frac13(k^2-2)\,\ellK^2(k)\,.
\end{equation}
\end{remark}

The nontrivial exercise in connection with this example is a
derivation of the analog to these rules for Heun's group $\G_s$ and
\eqref{Psi-s}. In other words, these rules, as a complete set of the
associated Picard--Fuchs equations, will be an `$s$-version' of the
following closed chain written in the `$\tau$-representation':
\begin{equation}\label{chain}
\dot s\dashrightarrow \vartheta,\theta\mbox{-forms}, \qquad \ddot
s\dashrightarrow \vartheta,\theta\mbox{-} \mbox{\ and\ }
\eta\mbox{-object},\qquad \dddot{\smash[b]{s}}
\dashrightarrow \bcQ(s)\,.
\end{equation}
Section~\ref{con} is devoted to explanations of this chain in a more
general context.

We close this section with a hyperelliptic example promised in
Sect.~\ref{Hyper}. From equations \eqref{tmp} and \eqref{sK} it
follows that
\begin{equation*}\label{temp}
\frac{16\,s}{s-1}\,\sqrt{x\,(x-1)\,}=\pm\sqrt{s\,(s-1)(s-9)\,}\,.
\end{equation*}
Remembering now that $s=z^2$, we see that this relation  leads us to
a hyperelliptic subcase of the family \eqref{Wab}, that is
\eqref{hyper}.

\begin{proposition}{\label{P5}
The hyperelliptic curve}
\begin{equation}\label{W2ab}
\begin{aligned}
v^2&=z^5-10\,z^3+9\,z\\&=z\,(z-1)(z+1)(z-3)(z+3)
\end{aligned}
\end{equation}
is described by Fuchsian equation \eqref{Qz} and has the following
parametrization\/$:$
$$
z=3\,\frac{\vartheta_3^2(3\,\tau)}{\vartheta_3^2(\tau)}\,,\qquad
v=48\sqrt{3\,}\ri
\frac{\vartheta_3^3(3\,\tau)}{\vartheta_3^3(\tau)}\,
\frac{\vartheta_2^2(\tau)\,\vartheta_4^2(\tau)}
{9\,\vartheta_3^4(3\,\tau)-\vartheta_3^4(\tau)}\,.
$$
\end{proposition}

\begin{example}\label{E56}
We can continue tower for the genus $g=2$ hyperelliptic curve
\eqref{W2ab} and find representation for  its uniformizing group
$\G$. General recipe to parabolic uniformization of hyperelliptic
Hauptmoduln is expounded in \cite[Sect.~6.2]{br}. A nontrivial
exercise (using Proposition~\ref{P5} and $\theta$-constant
techniques above)  is to construct the
$\vartheta,\theta$-representations for root functions $\sqrt{z\pm
1\,}$, $\sqrt{z\pm3\,}$ which are analogs of Weierstrassian
single-valued functions $\sqrt{\wp-e_k}=\sigma_k$. We note that
\eqref{W2ab} is readily simplified into the curve $v^2=z^6-1$.
\end{example}

\section{Hauptmoduln and analytic connections\label{con}}

The objectives to be pursued by Sect.~\ref{Fuchs} (especially around
formulae \eqref{Psi-s}--\eqref{chain}) were not only to exhibit way
of getting analytic formulae. All this material was substantially
intended for preparatory illustration to an important general
construction outlined in the chain \eqref{chain}. This is because
the number of independent differential consequences is not infinite;
they unify into a whole the basic aggregates of the theory:
Hauptmoduln $x,z,s,\...$, solutions $\Psi$ (possibly in terms of
$_2F_1$), and the primary differential $\Psi^2=\dot x$ as a weight-2
modular form. Closer motivation and explanations for that
interrelations are as follows.

\subsection{Differential structures on Riemann
surfaces\label{diffprop}}

Every Riemann surface $\R$ of finite analytic type (both genus and
number of punctures are finite \cite{bers,lehner}) is described by
the 2nd order equation \eqref{pqPsi}. Therefore we have the two
objects $\psi_1^{}(x)$, $\psi_2^{}(x)$ but their differential
closure requires, at most, the two more objects $\psi_1'(x)$,
$\psi_2'(x)$. A Wronskian relation sifts three of them as
differentially independent because
$$
\psi_1^{}\,\psi_2'-\psi_2^{}\,\psi_1'
=\exp\!\left(\![2]{-}\![4]\int\![4]p\,dx\right)
$$
is the known function. We may also think of the uniformizing
function $x=\chi(\tau)$ as a local coordinate transition
$x\rightleftarrows\tau$; of course, with a corresponding change of
complex structure preserving the $\R$ itself. Hence there should
exist the single-valued \mbox{`$\tau$-representations'} for
equivalents of these three base differential classes on $\R$ and the
$\R$ itself is considered now as a 1-complex-dimensional analytic
manifold. Clearly, we speak here of only analytic (not metric!)
objects since $\R$ is already assumed to be given.

Insomuch as all analytic tensor fields $T$ on Riemann surfaces
reduce to analytic \mbox{$k$-differentials} $T(\tau)\,d\tau^k$, we
can construct them from an object defining completely our $\R$, \ie,
the factor topology on $\overline{\Hp\!/\G}$. This is of course the
fundamental function $\bcQ(x,y)$ because it defines the above
mentioned $\psi$ and $\psi'$ in `$x$-representation'. It also
generates the scalar (automorphic function) $x(\tau)$ on this $\R$
through the $\bo{\mathfrak{D}}$-derivative $[x,\tau]=\bcQ(x,y)$; any
other equivalent Hauptmodul has the form $\frac{A\,x+B}{C\,x+D}$.
The simplest 1-tensor (call it $\g1$) is an analytic covariant and
general way of building such an object is to take a meromorphic
differential $\g1=R(x,y)\,dx$; our interest  now is only a
meromorphic analysis on $\R$. This differential, as an Abelian one,
is always the primary object $\dot x(\tau)$ (weight-2 automorphic
form) multiplied by any other scalar: $\g1(\tau)=R(x,y)\,\dot x$. In
the sections that follow we shall restrict our consideration to the
case $\bcQ=\bcQ(x)$ (genus $g=0$) and higher genera will be
considered elsewhere.

Motivated by the desire to build the differential apparatus on $\R$,
we need  to use derivatives of the Hauptmodul-scalar on $\R$ in
order to introduce the covariant differentiation, say, of
differentials $\nabla=\partial_u-\Gamma(u)$, with the standard
transformation law $\Gamma\mapsto\widetilde\Gamma$:
$$
\widetilde\Gamma(\p)\,d\tilde u=
\Gamma(\p)\,du-d\ln\frac{d\tilde u}{du}\,.
$$

There are many realizations of complex structures $u\mapsto \tilde
u$ for a given $\R$ but universal cover for all finite $\R$'s is
$\Hp$; this being so, we put $u\in\Hp$ and, renaming
$u\dashrightarrow \tau$, impose the projective structure
$\tau\mapsto\tilde\tau=\frac{\bo{a}\,\tau+\bo{b}}
{\bo{c}\,\tau+\bo{d}}$,
where $\big(\begin{smallmatrix}\bo{a}&\bo{b}\\
\bo{c}&\bo{d}
\end{smallmatrix}\big)\in\mathrm{PSL}_2(\mathbb{R})$. It follows
that we have to construct the coordinate representation
$\Gamma(\tau)$ of the connection object $\Gamma(\p)$ such that
$\Gamma(\tau)$---a function on $\Hp$---respects the factor topology
on $\overline{\HG}$ and transforms like a linear (affine) connection
does:
\begin{equation}\label{law}
\widetilde\Gamma(\p)=(\bo{c}\,\tau+\bo{d})^2\,
\Gamma(\p)+
2\,\bo{c}\,(\bo{c}\,\tau+\bo{d})\,,\qquad\p\in\R\,.
\end{equation}
On the other hand, we know that any scalar object $x(\p)=\tilde
x(\p)$ on arbitrary $\R$ is defined by the two structure properties:
\begin{equation*}
\tilde x\Big(\Mfrac{\bo{a}\,\tau+\bo{b}}
{\bo{c}\,\tau+\bo{d}} \Big)=
x(\tau)
\end{equation*}
and automorphism  of the functional representation $\chi(\tau)$ for
$x(\p)=\chi(\tau)$:
\begin{equation}\label{scalar}
\chi\mbig[7](\Mfrac{\alpha\,\tau+\beta}
{\gamma\,\tau+\delta} \mbig[7])=\chi(\tau)\quad \forall\:
\mbig[7](\begin{matrix}\alpha&\!\!\beta\\
\gamma&\!\!\delta \end{matrix}\,\mbig[7])\in\G_x
\end{equation}
(we supposed that $\pi_{\!1}(\R)=\G_x$). Hence it follows that the
$\Gamma$-object can be built with the help of differential or scalar
by the rule
\begin{equation}\label{dxt}
\Gamma(\tau)=\frac{d}{d\tau}\!\ln \g1=
\frac{d}{d\tau}\!\ln R(x)\,\dot x(\tau)\,.
\end{equation}
To put it differently, we are interested in $\Gamma$'s compatible
with scalars  on $\R$.

The above-mentioned differential closure is already available
because the Schwarz object $\{x,\tau\}\,d\tau^2$ becomes (as it
follows from \eqref{tensor}) a tensor when complex structure is a
projective one. Hence $\{x,\tau\}\,d\tau^2=\bcQ(x)\,dx^2$ and,
finally, the transition between $x$- and
\mbox{$\tau$-representations} can be represented in the form of the
following  equivalence
$$
\mbig\{\psi_1^{}(x)\,,\:\psi_2^{}(x)\,,\:\psi_1'(x)
\mbig\}\quad\hhence\quad
\mbig\{\chi(\tau)\,,\:\dot \chi(\tau)\,,\:\Gamma(\tau) \mbig\}\,,
$$
where
\begin{equation}\label{Gt}
\tau=\frac{\psi_2^{}(x)}{\psi_1^{}(x)}\,,\qquad\dot
\chi(\tau)=\bE\,\psi_1^2(x)\,,\qquad\Gamma(\tau)=
2\,\bE\,\psi_1^{}\psi_1'+(p+{\ln'}\!R)\,\bE\,\psi_1^2
\end{equation}
and $\bE\DEF\exp\![3]\int\![3] p\,dx$; if the normal form
\eqref{fuchs} has been chosen then one puts here $(p,\bE)=(0,1)$. We
thus conclude finally.
\begin{itemize}
\item Complete data set for the analytic theory on $\R$ is
    defined by the relations \eqref{Gt} and is controlled
    \emph{not by one} classical equation \eqref{pqPsi}, but by
    the two fundamental \odes\ for the $\psi$-function and its
    derivative $\psi'\FED\phi$:
\begin{equation}\label{base}
\begin{array}{c}
\ds\psi''+p\,\psi'+q\,\psi=0\,,\\[2ex]
\ds\phi''+(p-{\ln'}q)\,\phi'+(q+p'-p\,{\ln'}q)\,\phi=0\,.
\end{array}
\end{equation}
Choice of functions $R(x)$, under this setting, calibrates all
the $\Gamma$'s.
\end{itemize}
The second of Eqs.~\eqref{base} is also of Fuchsian class and it is
a very interesting question on role of its monodromy in the theory.
This point was not arisen in the literature even for the classical
case $\bo\Gamma(2)$ for which one already has an exhaustive
collection of formulae \eqref{KE1}--\eqref{KE4}. For example, under
normalization $\phi_1=\Psi_1'$, $\phi_2=(\tau\,\Psi_1)'$ it is not
difficult to see that matrix monodromy $\G^{(\phi)}$ of the
$\phi$-equation is determined by $\tau$-transformation of the
single-valued function
$$
\frac{\phi_2(\tau)}{\phi_1(\tau)}=
2\,\frac{\dot \chi(\tau)}{\ddot \chi(\tau)}+
\tau=\cdots
$$
which is the same as
$$
\cdots=\frac{2}{\Gamma(\tau)}+\tau
$$
under the choice $R(x)=1$ in \eqref{dxt}. Self-suggested questions
here are the relation of the monodromy $\G^{(\phi)}$ with $\G_x$,
its genus (always not finite?), and more precise/functional
relationship with $\Gamma(\tau)$ and automorphic function
$x=\chi(\tau)$.

Of course, we may (for convenience sake) freely change the basis
$\{\psi_1^{},\psi_2^{},\phi_1^{},\phi_2^{}\}$ over $\mathbb{C}$ or
$\mathbb{C}(x)$ and rewrite (in an equivalent way) the theory in
terms of the four 1st order \odes. For example, formulae
\eqref{KE1}--\eqref{KE4} and \eqref{var} illustrate this
construction when group is $\bo\Gamma(2)$.

\begin{remark}\label{R8}
The geometric treatment to the famous Chazy equation
$\pi\,\dddot{\smash[b]{\eta}}=12\,\ri(2\,\eta\,\ddot\eta-
3\,\dot\eta^2)$ as one defining the \mbox{$\Gamma$-connection} for
group $\bo\Gamma(1)$ was given by Dubrovin in lectures
\cite{dubrovin}. In the same place he showed that other equations
found by Chazy fit also in such a scheme. See Appendix C of
\cite{dubrovin} for more examples and voluminous references. The
following simple example shows that Dubrovin's
\mbox{$\bo\Gamma(1)$-connection} is compatible with scalar.
\end{remark}

\begin{example}\label{E10}
In the case of group $\bo\Gamma(1)$ we may take $x=J(\tau)$. The
simplest idea for a basic weight-2 automorphic form $\g1$ is to take
the ratio of modular Weierstrass' forms $\g1=\frac{g_3^{}}{g_2^{}}$.
This form does not contradict to usage of the
$\bo\Gamma(1)$-Hauptmodul and we may put
\begin{equation}\label{gG}
J(\tau)\,,\qquad \g1(\tau)\DEF\frac{\pi\ri}{36}\,\frac{\dot
J}{J}\,,\qquad
\Gamma(\tau)\DEF\frac{\ddot J}{\dot J}\,.
\end{equation}
As a references source we give  relations to the standard functions,
\ie, representations for standard \emph{holomorphic} forms on
$\bo\Gamma(1)$ through the objects \eqref{gG}:
$$
\g2(\tau)=\frac{27\,J}{J-1}\,\g1^2\,,\qquad
\g3(\tau)=\frac{27\,J}{J-1}\,\g1^3\,,
\qquad
\eta(\tau)=\frac{\pi}{4\ri}\,\Gamma+\frac32\frac{7\,J-4}{J-1}\,\g1
$$
and
$$
\nabla
\g1=(\partial_\tau-\Gamma)\,\g1=\frac{6}{\pi\ri}\,
\frac{J+2}{J-1}\,\g1^2\,.
$$
The connection $\Gamma$ is not holomorphic everywhere in $\Hp$ but
all the $\Gamma$'s are defined up to a 1-differential.
Subtracting/adding a certain differential we obtain the holomorphic
connection
\begin{equation}\label{J1}
\Gamma\dashrightarrow\Gamma-\frac16\,\frac{d}{d\tau}\!\ln J^4(J-1)^3=
\frac{4\ri}{\pi}\eta(\tau)\,.
\end{equation}
\end{example}

\begin{remark}[\sf Example]\label{R9}
As a continuation of this example we may  write down the
$\bo\Gamma(1)$-equations \eqref{base} and complete set of their
solutions in terms of special functions (not merely in terms of
${}_2F_1$-series). Curiously, the group $\bo\Gamma(1)$ is not a
lesser known group than $\bo\Gamma(2)$ but no formulae solving this
problem  have appeared in the literature; we mean direct/inverse
modular transitions between $\{\psi(x),\phi(x)\}$ and $\{\tau\}$  as
analogs of \eqref{KE1}--\eqref{KE4}. Partially these points  were
solved in \cite{br4}.
\end{remark}

\subsection{\odes\ satisfied by analytic connections on
$\R$'s\label{anal}}

Since the general classes of Dubrovin's $\Gamma$-equations
\cite[p.~152, {\bf Exercise C.3}]{dubrovin} respect the most general
projective structure $\mathrm{GL}_2(\mathbb{C})$, their discrete
symmetries (if any) must not necessarily be Fuchsian. Without any
specification these may  formally be non-finitely generated groups
or described even by not Fuchsian equations\footnote{A
counterexample: formula (\emph{C}.44) in {\bf Exercise C.4} on
p.~155 in \cite{dubrovin}. Connections and Hauptmoduln associated to
that equations are not bound to be single-valued functions;  factor
topologies, automorphic identities, and groups are not irrelevant in
such a kind of examples.} (to say nothing of monodromies of Fuchsian
type). In other words, projectively correct $\Gamma$-equations do
not touch on the question of finiteness of genus and compatibility
of connections with Hauptmoduln-scalars, whereas any finite $\R$ is
associated, as described above, with such $\Gamma$'s (hence
Dubrovin's $\Gamma$'s as well) and we have actually had now a large
number of examples with their monodromies known to be Fuchsian.

\begin{theorem}\label{T10}
Let $\bo{\mathfrak{T}}$ be a genus zero uniformizing orbifold
defined by the Fuchsian equation $[x,\tau]=\bcQ(x)$ having only
parabolic singularities $x=E_k$. Then the expression
\begin{equation}\label{Gamma}
\Gamma=\frac{d}{d\tau}\!\ln\dot x
\end{equation}
is an analytic connection $($respecting the factor topology on
$\overline{\HG}$\/$)$ everywhere holomorphic on $\Hp
\cup\{\infty\}$. It satisfies a 3rd order polynomial \ode \
$\Xi(\dddot{\smash[b]{\,\Gamma}},\ddot{\mathrm\Gamma},
\dot{\mathrm\Gamma},\Gamma)=0$ with constant coefficients. Any other
analytic connection also satisfies an \ode\ of such a kind. General
solutions to these equations are constructed by the scheme
\begin{equation}\label{general}
\Gamma(\tau)\quad\mapsto\quad\frac{\Gamma
\big(\frac{\bo{a}\,\tau+\bo{b}}
{\bo{c}\,\tau+\bo{d}}\big)}
{(\bo{c}\,\tau+\bo{d})^2}-\frac{2\,\bo{c}}
{\bo{c}\,\tau+\bo{d}}\,.
\end{equation}
\end{theorem}

\begin{pf}
The statement about factor topology is obvious from the property
\eqref{scalar} of $x=\chi(\tau)$ of being a scalar. Since
$\bo{\mathfrak{T}}$ is an orbifold of finite genus, the existence
domain of $\chi(\tau)$ is an interior of a circle or a half-plain.
Normalize this domain to be $\Hp$. For all $\tau,\tau_0^{}\in\Hp$ we
have a convergent series representation
$\chi(\tau)=x_0^{}+a\,(\tau-\tau_0^{})+\cdots$, where $a\ne0$. It
follows that \eqref{Gamma} is holomorphic everywhere in $\Hp$. For
the infinite point we make the standard change $\tau\mapsto q$ of
the local parameter $q=\re^{\pi\ri\tau}$, where $\tau\to\ri\infty$.
Then we may write $x=E+a\,q^n+\cdots$, where $E$, $a$, and
$n\in\mathbb{Z}$ depend on the local monodromy $\G_x$. Taking into
account that $dq=\pi\ri q\,d\tau$, we derive that \eqref{Gamma} is
again holomorphic as $q\to0$. Zeroes and behavior of connection on
the real axis are not well defined since its transformation law is
inhomogeneous.

Let us denote $\Omega\DEF\dot{\mathrm{\Gamma}}-\frac12\,\Gamma^2$.
The equation $[x,\tau]=\bcQ(x)$ determining the scalar $x$ becomes
the tensor identity $\bcQ\cdot\dot x^2=\Omega$. Applying  $\nabla$,
we get the two equations:
\begin{equation}\label{nabla}
\Omega=\bcQ\cdot(\nabla x)^2\,,\qquad
\nabla\Omega=\nabla\bcQ\cdot(\nabla x)^2= \bcQ'\!
\cdot(\nabla x)^3\,,
\end{equation}
since $\nabla x=\dot x$ and $\nabla\dot x\equiv0$. Elimination of
$\nabla x$ gives the identity
\begin{equation}\label{23}
\frac{(\nabla\Omega)^2}{\Omega^3}=\frac{\bcQ'^2}{\bcQ^3}\,.
\end{equation}
The second covariant derivative $\nabla^2\Omega=\bcQ''\!
\cdot(\nabla
x)^4$ yields yet another identity:
\begin{equation}\label{''}
\frac{\nabla^2\Omega}{\Omega^2}=\frac{\bcQ''}{\bcQ^2}
\end{equation}
and $\nabla$-derivatives are understood here to be equal to
\begin{equation}\label{GO}
\nabla\Omega\DEF\dot{\mathrm{\Omega}}-2\,\Gamma\,\Omega\,,\qquad
\nabla^2\Omega\DEF(\partial_\tau-3\,\Gamma)
(\partial_\tau-2\,\Gamma)\,\Omega\,.
\end{equation}
Since identities \eqref{nabla}, \eqref{23} are of invariant (scalar)
type, the sought-for $\Gamma$-equation results from elimination of
$x$. This will be a 3rd order \ode\
$\Xi(\Omega,\nabla\Omega,\nabla^2\Omega)=0$ with constant
coefficients. Notice incidentally that such an equation does still
exist if group $\G_x$ is not a 1st kind Fuchsian one, \ie, when
$\bcQ(x)$ determines a Schottky group or even has non-correct
(`bad') accessory parameters.

Let $\gamma$ be any other connection. Then $\Gamma=\gamma+R(x)\,\dot
x$, where $R(x)\,\dot x$ is a certain \mbox{1-differential}. For
example, if original Fuchsian equation has non-canonical form
\eqref{pqPsi}, it is convenient to put $R(x)=p(x)$, that is
$\Gamma=\gamma+p(x)\,\dot x$. Hence we redefine $\Omega$ as
$$
\Omega=\partial_\tau(\gamma+R\,\dot x)-
\frac12\,(\gamma+R\,\dot x)^2=
\mbig[7](\dot\gamma-\frac12\,\gamma^2\mbig[7])+
\mbig[7](R'+\frac12\,R^2\mbig[7])\dot x^2\,,
$$
since $\ddot x=(\gamma+R\,\dot x)\,\dot x$. Denoting
$\omega\DEF\dot\gamma-\frac12\,\gamma^2$ one gets
\begin{equation}\label{222}
\omega=\mbig[7](\bcQ-R'-\frac12\,R^2\mbig[7])\cdot\dot x^2
\qquad\hence\qquad\omega=Q\cdot(\widetilde\nabla x)^2\,,
\end{equation}
where $Q\DEF\bcQ-R'-\frac12\,R^2$ and $\widetilde\nabla$ signifies
the differentiation by means of connection $\gamma$. Clearly, $x$ is
a `flat' coordinate only with respect to the `old' connection
\eqref{Gamma}, \ie, $\widetilde\nabla^2  x\not\equiv0$ now, and  we
have $\widetilde\nabla^2 x=R\cdot(\widetilde\nabla x)^2$ instead of
$\nabla^2x\equiv0$. Therefore $\widetilde\nabla$-derivatives of
\eqref{222} give the two identities
$$
\widetilde\nabla\omega=\big(Q'+2\,R\,Q\big)\,
(\widetilde\nabla x)^3\,,
\qquad
\widetilde\nabla^2\omega=\big(
Q''+5\,R\,Q'+2\,R'Q+6\,R^2Q\big)\,(\widetilde\nabla x)^4
$$
which are the generalizations of Eqs.~\eqref{23}--\eqref{''}:
$$
\frac{(\widetilde\nabla\omega)^2}{\omega^3}=
\frac{(Q'+2\,R\,Q)^2}{Q^3}\,,
\qquad
\frac{\widetilde\nabla^2\omega}{\omega^2}=
\frac{1}{Q^2}\,\big(Q''+5\,R\,Q'+2\,R'Q+6\,R^2Q\big)\,.
$$
As before, the equation
$\Xi(\dddot{\smash[b]{\gamma}},\ddot\gamma,\dot\gamma,\gamma)=0$
follows by elimination of $x$ and
\begin{equation}\label{gamma}
\widetilde\nabla\omega\DEF(\partial_\tau
-2\,\gamma)\!\left(\dot\gamma-\frac12\,\gamma^2\right),\qquad
\widetilde\nabla^2\omega\DEF(\partial_\tau
-3\,\gamma)(\partial_\tau
-2\,\gamma)\!\left(\dot\gamma-\frac12\,\gamma^2\right).
\end{equation}
Transformation law \eqref{law} for all the $\Gamma$'s entails the
formula \eqref{general} for solutions to these equations wherein
$\Gamma(\tau)$ is any particular solution, \eg, \eqref{Gamma}, if
Hauptmodul $x=\chi(\tau)$ has been given.\hfill $\blacksquare$
\end{pf}

Notice that in the canonical case $(p,R,\bE)=(0,1,1)$ the
logarithmic derivative $2\ln_\tau\!\Psi$ satisfies the same equation
as $\Gamma$ does; according to the 3rd formula in \eqref{Gt}, we
have
$$
\Gamma(\tau)=2\,\Psi\,\Psi'=2\,\frac{\Psi_\tau}{\Psi}
$$
(this is yet another motivation for introducing the
$\Psi'$-derivative). It should be noted here that the 3rd order
$\gamma$-\ode\ may turn out to be simpler if the $\gamma$-definition
corresponds not to $\G_x$ but to a wider group. For example, the
connection
$$
\gamma=\frac{d}{d\tau}\!\ln\dot s -
\frac32\,\frac{s-5}{s\,(s-9)}\,\dot
s\,,
$$
has been constructed formally as one for the Heun group $\G_s$
\eqref{star} with an addition of the differential $R(s)\dot s$. But
this $\gamma$ is in fact a hidden form of the
$\bo\Gamma(1)$-connection satisfying the equation
$\dddot{\smash[b]{\gamma}}= 6\,\gamma\,\ddot\gamma-9\,\dot\gamma^2$;
proof is a calculation with use of \eqref{J} and \eqref{J1}. In
other words
\begin{itemize}
\item \emph{Connections, along with uniformizing Hauptmoduln
    $x=\chi(\tau)$, also form towers according to a tower of
    subgroups}.
\end{itemize}

\begin{example}\label{E11}
The connection for Legendre's modulus
$k(\tau)=\frac{\vartheta_2^2(\tau)}{\vartheta_3^2(\tau)}$ with
monodromy $\G_u$ defined by equation \eqref{lemn}.  From
Theorem~\ref{T10} we obtain the following $\Gamma$-equation:
$$
\begin{aligned}
&A^8-8\,\Omega\,(B-352\,\Omega^2)\, A^6+24\,\Omega^2\,
(B^2-260\,\Omega^2\,B-368\,\Omega^4)\,A^4\\
&\phantom{A^8}-32\,\Omega^3\,(B^3-129\,\Omega^2\,B^2-168\,\Omega
^4\,B+944\,\Omega^6)\,A^2+16\,\Omega^4\,
(B^2-20\,\Omega^2\,B-80\,\Omega^4)^2=0\,,
\end{aligned}
$$
where $A\DEF\nabla\Omega$, $B\DEF\nabla^2\Omega$, and $\Omega$ as in
\eqref{GO} and
$$
\Gamma
=\frac{d}{d\tau}\!\ln \!\frac{d}{d\tau}k(\tau)
=\frac{4\ri}{\pi}\,\eta+
\frac16\,\pi\ri\big(\vartheta_3^4-5\,\vartheta_2^4\big)\,.
$$
Renormalization of this $\Gamma$ into connection
$$
\begin{aligned}
\Gamma\dashrightarrow\gamma=\Gamma-\frac{d}{d\tau}\!\ln k\,(k^2-1)&=
\frac{d}{d\tau}\!\ln\vartheta_2^4\\
&=\frac{4\ri}{\pi}\,\eta+
\frac13\,\pi\ri\big(2\,\vartheta_2^4-\vartheta_3^4\big)\,,
\end{aligned}
$$
which corresponds to original Chudnovsky equation \eqref{I}, yields
more compact equation:
$$
(2\,\dot\gamma-\gamma^2)\,\dddot{\smash[b]{\gamma}}
=2\,\ddot\gamma\,(\ddot\gamma-\gamma^3)-
\dot\gamma^2(2\,\dot\gamma-3\,\gamma^2)\,.
$$
\end{example}

\begin{example}\label{E12}
The $\gamma$-equation for  Chudnovsky equation \eqref{II}. It admits
the compact form
$$
2\,\omega\,\big(\widetilde\nabla^2\omega+6\,\omega^2\big)^2
=\big(2\,\widetilde\nabla^2\omega+15\,\omega^2\big)\,
(\widetilde\nabla\omega)^2
\,,
$$
where $\widetilde\nabla\omega$, $\widetilde\nabla^2\omega$ are
determined by \eqref{gamma} and the scalar compatible $\gamma$ is
defined here as
$$
\gamma=\frac{d}{d\tau}\!\ln\dot x-\frac{d}{d\tau}\!\ln
\big\{(x+1)^3-1\big\}\,.
$$
\end{example}

We do not display here formulae related to most interesting equation
\eqref{heun} since its $\Gamma,\gamma$-equations have somewhat
cumbrous form (it is a good exercise to list these \odes\ for all
Chudnovsky's equations). It is easily derivable but we note in
passing that these examples, \ie, $\G_u$ and $\G_s$, show that  the
direct search for $\Gamma$-equations in form of  Dubrovin's tensor
$\{\Omega$, $\nabla\Omega$, $\nabla^2\Omega,\,\...\}$-Ans\"atze
\cite{dubrovin} can be very difficult problem, whereas all these
equations will known to be consequences of those coming from the
scalars; all of them are at our disposal once we have had a
Hauptmodul(n). Moreover, normalization of connections by a
differential $R(x)\dot x$ significantly affects the size of
equation. For example, the natural but `unlucky'
$\bo\Gamma(1)$-definition \eqref{gG} instead of \eqref{J1} leads to
useless equations like
$$
\begin{aligned} 2642368542992676\,B^6-
264634415204382300\,B^5\,\Omega^2+\cdots&\\
{}+26676039242547325440000\,\Omega^{12}&=0
\end{aligned}
$$
with the same meaning for $A$, $B$, and $\Omega$ as above.

If Hauptmodul has elliptic singularities then construction of the
holomorphic connection is performed by a simple reproducing of what
we have done in {\sf Example~\ref{E10}}, that is \eqref{J1}.

Let the zero genus Hauptmodul $x(\tau)$ have an $N$-order conical
point $x=C$; that is $x=C+a\,(\tau-\tau_0^{})^N+\cdots$. In order to
compensate singularities coming from the repeated points of
$x(\tau)$ we  add to connection \eqref{Gamma} the differential
$\frac{1-N}{N}\frac{\dot x}{x-C}$ for each such point. The $\Gamma$
remains regular at infinity.

\begin{proposition}\label{P6}
Everywhere in $\Hp \cup\{\infty\}$ holomorphic affine connection for
a zero genus orbifold defined by the Fuchsian equation
$[x,\tau]=\bcQ(x)$ with $N_k$-order elliptic singularities at points
$x=C_k$ is determined by the expression\/$:$
$$
\Gamma=\frac{d}{d\tau}\!\ln\dot
x-\sum_k\frac{N_k-1}{N_k}\,\frac{d}{d\tau}\!\ln(x-C_k)\,.
$$
\end{proposition}

Although connections are not uniquely defined we can partially
remove the ambiguity by normalizing location of the three points at
fixed cusps, say, at  $x=\{0,1,\infty\}$.

Study of nonlinear differential equations associated with certain
modular forms and generalizations of Chazy's equations were
initiated by M.~Ablowitz et all in the nineties in connection with
mathematical physics problems including magnetic monopoles
\cite{hitchin2}, self-dual Yang--Mills and Einstein equations
\cite{hitchin}, as well as topological field theories
\cite{dubrovin}. Recently they again attracted attention
\cite{ablowitz} and  the nice work \cite{maier} by R.~Maier
provides explicit examples related to some low level groups
$\Gamma_0(N)$, hypergeometric equations, and number-theoretic
treatments. The works \cite{ablowitz,maier} provide also voluminous
references along these lines.

We left  some remarks and examples in this work as exercises because
the stream of consequences, including hyperelliptic, may be
increased considerably. The abundance of  towers, Hauptmoduln,
connections, \odes, and groups requires their further
classification, and with it the unification of ways of getting the
$\theta$-formulae and groups.

\section{Remarks on Abelian integrals and equations on
tori\label{final}}

Described differential `$\theta$-machinery' makes it possible to
include immediately into analysis Abelian integrals: holomorphic,
meromorphic, and logarithmic. Indeed, closedness of the construction
\eqref{Gt}--\eqref{base} and differential calculus proposed in
Sect.~\ref{diffprop} are actually just a `semi-true` because, to be
fully consistent, we should also extend Fuchsian theory into these
types of integrals. This goes far beyond the scope of the present
work and will be the subject matter of a separate study.

Insomuch as many of the curves we have considered (probably all)
admit representation in form of covers over tori, we can derive a
large family of explicitly solvable Fuchsian equations on tori; in
doing so fundamental logarithmic Abelian integral and all the
meromorphic objects on our Riemann surfaces can be described in
terms of Jacobi's $\theta$-functions and constants; a great
effectivization of the theory. As we have already mentioned in
Introduction, no such representations are presently available % and
%direct statement of this problem, to the best of our knowledge,
%has not been  even raised in the literature
in spite of numerous examples of modular curves and the well-studied
subgroups of $\bo\Gamma(1)$. We  touched briefly on this problem in
\cite{br}, considered there one special example of a meromorphic
integral, and exhibited in \cite[Proposition~10]{br5} the first
example of analytic formula for everywhere holomorphic and
additively automorphic object on $\mathbb{H}^+$. Abelian integrals,
Fuchsian equations on tori, and higher genus integer $q$-series
certainly merit thorough investigation and will be fully considered
in a continuation of this work. Anticipating such a kind results by
explicit formulae, we consider briefly the hyperelliptic family
\eqref{hyper}, in particular, curve \eqref{W2ab}.

It is common knowledge that all such curves admit a reduction of
their holomorphic integrals  into the elliptic ones. Most known
reduction formula is due to Jacobi, concerns with curve
\eqref{W2ab}, and has in our case the form
\begin{equation}\label{pzk}
\wp(\bo{u})=\frac{(z\pm\varkappa)^2}{(z\mp
\varkappa)^2}
\qquad (\varkappa\DEF\ri\sqrt{3})\,;
\end{equation}
it respects also the relation $\wp'(\bo{u})^2=4\,\wp^3(\bo{u})-4$
when $z(\tau)$ satisfies \eqref{Qz}. It follows (nontrivial usage of
Lemma~\ref{L1}) that the two nice formulae describe the theory of
this curve:
$$
\sqrt{\mp\varkappa}\cdot
d\bo{u}=\frac{(z\pm\varkappa)\,dz}{\sqrt{z^5-10\,z^3+9\,z}}\,,
\qquad
[\bo{u},\tau]=-6\,\frac{2\,\wp^3(\bo{u})+1}
{\wp^2(\bo{u})\,\wp'(\bo{u})^2}
$$
(computation details will be presented elsewhere). Let us specify
what does this mean. Both the holomorphic integrals
$\bo{u}^\pm=\bo{u}(\tau)$ defined by the equations above (Fuchsian
equations on a torus) are the globally single-valued analytic
functions on $\Hp$ because equation \eqref{pzk} is obviously none
other than an equivalent of the  Riemann surface defined by
hyperelliptic form \eqref{W2ab}.

Yet another example is related to the Heun--Ap\'ery Eq.~\eqref{heun}
and we sketch the broad outlines of the theory (taken from
\cite{br5}):
\begin{equation}\label{U}
[\mathfrak{u},\tau]=-2\,\wp(2\mathfrak{u})-\frac83
\qquad\hhence\qquad
[\bo{\mathfrak{u}},\tau]=
-2\,\wp(2\bo{\mathfrak{u}}|\varepsilon)-
\frac13\,\pi^2\vartheta_2^4(\varepsilon),
\end{equation}
where $\mathfrak{u}=\omega\,\bo{\mathfrak{u}}$ and Weierstrass'
$\wp(\mathfrak{u})=\wp(\mathfrak{u};a,b)=
\wp(\mathfrak{u}|\omega,\varepsilon\,\omega)$ corresponds to
invariants $(a,b)=\mbig(\frac{292}{3},\frac{4760}{27} \mbig)$ and
therefore \mbox{$J(\varepsilon)=\frac{73^3}{2^4 3^7}$}
\cite[p.~185]{chud2}. We failed to find out  in reference books on
imaginary quadratic number fields (\eg, \cite{weber}, \cite{watson})
an exact value for $\varepsilon\approx
1.5634019226921973634612986241\cdot\ri$, so this (Chudnovsky's)
torus is perhaps the very exceptional indeed; see \cite{br5}.
Monodromy of Eq.~\eqref{U} equation is generated by the two
transformations $\mathfrak{a}(\bo{\mathfrak{u}})=
\bo{\mathfrak{u}}+1$ and $\mathfrak{b}(\bo{\mathfrak{u}})=
\bo{\mathfrak{u}}+\varepsilon$; the loop surrounding
$\bo{\mathfrak{u}}=0$ is equivalent to the transformation
$\mathfrak{a}\,\mathfrak{b}\,\mathfrak{a}^{\sm1}\,
\mathfrak{b}^{\sm1}$. A result of the theory is that equations
\eqref{U} and \eqref{heun} are turnable one to another by the very
simple substitution
\begin{equation}\label{subs}
s=\wp(\mathfrak{u})+\frac{10}{3}
\end{equation}
and therefore the global monodromy $\G_\mathfrak{u}$ is an index two
2-generated subgroup of Heun--Ap\'ery's one $\G_s$ described in
Theorem~\ref{T9}. By construction,  solution to Eq.~\eqref{U} should
be sought-for as an additively automorphic and globally
single-valued analytic function on $\Hp$. Partially, an explicit
representation to $\mathfrak{u}(\tau)$ is in fact the latter formula
$$
\mathfrak{u}(\tau)=\wp^{\sm1}\!\bigg(9\frac{\vartheta_3^4(3\tau)}
{\vartheta_3^4(\tau)}-\frac{10}{3}\bigg).
$$
Remembering now that $s=z^2$, we may write an equivalent of
substitution \eqref{subs}:
\begin{equation}\label{unit}
z^2=\wp(\mathfrak{u})+\frac{10}{3}
\end{equation}
which in turn coincides with a simplest  example of a nontrivial
($g>0$) 2-sheeted cover over tori, \ie, representations of a special
but very wide class of $\R$'s. It has the form
$z^2=\wp(\mathfrak{u})-e$, where $e$ is any of the standard
Weierstrass' branch points. Since
$\wp(\varepsilon\omega)=-\frac{10}{3}$, the cover \eqref{unit} is of
unit genus, has invariant $J=\frac{2197}{972}$ which being, as it
should, equal to $J$-invariant of the Dubrovin--Hitchin torus
\eqref{tor-yu} from whence Eq.~\eqref{U} has come in our
constructions.

Further remarkable consequence of arisen `toroidal covers' is a
corollary on the equivalence of Chudnovsky's equations, say
\eqref{I} and \eqref{III}. It follows at once that representation of
these two orbifolds through  punctured tori will yield a nontrivial
representation for the mutually transcendental covers
\mbox{$\Xi(\mathfrak{p},\mathfrak{u})=0$} of the punctured torus
\eqref{U} and a torus defined by Chudnovsky's equation \eqref{I}.
The latter has long been known \cite{keen,chud2} and has the form of
a `punctured lemniscate'
$[\mathfrak{p},\tau]=-2\wp(2\mathfrak{p}|\ri)$. It is transformed
into Chudnovsky's equation \eqref{I} by the further simple
substitution $x=\wp(\mathfrak{p}|\ri)$ (good exercise is to check
this). Correlating this substitution, \eqref{subs}, and \eqref{sx},
one can derive the sought-for function~%
$\Xi(\mathfrak{p},\mathfrak{u})$ (a simpler version than that from
\cite{br5}):
\begin{equation}\label{pu}
\mathfrak{u}\rightleftarrows\mathfrak{p}:\qquad
\omega\,\frac{\mbig[3](\wp(\mathfrak{u};a,b)+\frac13\mbig[3])^2-12}
{2\,\wp^2(\mathfrak{p};4,0)-1}=
\sqrt{8\,}\,\pi\,\vartheta_2\!\!\left(\Mfrac\varepsilon2
\right)\cdot\frac{\theta_4}{\theta_1}
\mbig[7](\frac{\mathfrak{u}}{2\omega}\Big|\varepsilon\mbig[7])
\end{equation}
($\omega=\frac12\pi\ri\vartheta_2^2(\varepsilon)\approx
0.539\,128\,911\,874\,910\,808\,859\,668\,749\!\...$). To put it
differently, as cannot well be imagined, the two `simple' Fuchsian
equations
$$
[\mathfrak{u},\tau]=-2\,\wp(2\mathfrak{u};a,b)-\frac83\,,\qquad
[\mathfrak{p},\tau]=-2\,\wp(2\mathfrak{p};4,0)
$$
are transformable into each other by the  transcendental and highly
non-obvious cover \eqref{pu}. In this regard they are also
integrable along with the four Chudnovsky's equations. The cover
\eqref{pu} is of course completely uniformized by single-valued
functions $\mathfrak{u}(\tau)$ and $\mathfrak{p}(\tau)$. It is a
very nontrivial exercise to investigate it directly, \ie, to
describe its branch schemes, Puiseux series, and compute its
(nontrivial) genus.

What we have done now is in effect the draft recipes of getting the
large number of analytic formulae including the very nontrivial
`toroidal towers'. But all this is just a corollary of simplest
PH-curves. The \Psix-equation \eqref{p6} is thus not the mere  rich
source of algebraic solutions. It generates infinite families of
explicitly uniformizable curves, their toroidal cover
representations, if any, and complete differential apparatus on
these orbifolds.

\end{document}